\documentclass[hidelinks,onefignum,onetabnum]{siamart250211}

\usepackage{lipsum}
\usepackage{amsfonts}
\usepackage{graphicx}
\usepackage{epstopdf}
\usepackage{algorithmic}
\ifpdf
  \DeclareGraphicsExtensions{.eps,.pdf,.png,.jpg}
\else
  \DeclareGraphicsExtensions{.eps}
\fi

\newsiamremark{remark}{Remark}
\newsiamremark{hypothesis}{Hypothesis}
\crefname{hypothesis}{Hypothesis}{Hypotheses}
\newsiamthm{claim}{Claim}
\def\Frac#1#2{\frac{\displaystyle{#1}}{\displaystyle{#2}}}
\def\dsp{\displaystyle}
\def\wt{\widetilde}

\def\Ai{{{\rm Ai}}}

\def\wt{\widetilde}
\def\tfrac#1#2{{{\lower.6ex
\hbox{$\scriptstyle#1$}}\over 
{\raise.7ex
\hbox{$\scriptstyle#2$}}}}

\def\tfrac#1#2{{{\lower.6ex
\hbox{$\scriptstyle#1$}}\over 
{\raise.7ex
\hbox{$\scriptstyle#2$}}}}

\headers{Computation of classical Gaussian quadratures}{A. Gil, J. Segura and N. M. Temme}

\title{Fast and accurate computation of classical Gaussian quadratures   
\thanks{Submitted to SIAM J. Sci. Comput.
\funding{The authors acknowledge support from Ministerio de Ciencia e Innovaci\'on, projects
 PID2021-127252NB-I00 (MCIN/AEI/10.13039/501100011033/ FEDER, UE) and
 PID2024-159583NB-I00 (MICIU/ AEI / 10.13039/501100011033 / FEDER, UE) 
 }}
 }

\author{
 Amparo Gil\thanks{Depto. de Matemática Aplicada y Ciencias de la Computaci\'{o}n, Universidad de Cantabria, Avda. de los Castros, s/n, 39005 Santander, Spain
	(\email{amparo.gil@unican.es}).}
\and Javier Segura\thanks{Depto. de Matemáticas, Estadística y Computaci\'{o}n, Universidad de Cantabria, Avda. de los Castros, s/n, 39005 Santander, Spain
	(\email{javier.segura@unican.es}).}
		\and N.M. Temme\thanks{Valkenierstraat 25, 1825 BD Alkmaar, the Netherlands (\email{nic@temme.net}).}
}

\usepackage{amsopn}

\ifpdf
\hypersetup{
  pdftitle={Fast and accurate computation of classical Gaussian quadratures},
  pdfauthor={A. Gil, J. Segura and N. M. Temme}
}
\fi

\usepackage{xcolor}
\usepackage{graphicx}
\definecolor{codepurple}{rgb}{0.58,0,0.82}
\usepackage{listings}
\definecolor{codegreen}{rgb}{28,172,0}
\usepackage[utf8]{inputenc}

\newcommand{\listtt}[1]{\lstinline|#1|}

\usepackage{dirtree}
\usepackage{framed}
\renewenvironment{framed}[1][\hsize]
  {\MakeFramed{\hsize#1\advance\hsize-\width \FrameRestore}}%
  {\endMakeFramed}

\begin{document}

\maketitle

\begin{abstract}
Algorithms for computing the classical Gaussian quadrature rules (Gauss--Jacobi, Gauss--Laguerre, and Gauss--Hermite) are presented, based on globally convergent fourth-order iterative methods combined with asymptotic approximations, which are applied in complementary regions of the parameter space. This approach yields methods that improve upon existing algorithms in speed, accuracy, and computational range. 
The MATLAB algorithm for Gauss--Jacobi is faster than previous methods and lifts the upper 
restrictions on the parameters imposed by those methods ($\alpha,\beta\le 5$); 
for example, for degrees up to $10^6$ all nodes and weights can be 
computed within the underflow limit 
for $-1<\alpha,\beta\le 30$, and the computable range of parameters is much larger for smaller degrees, limited only by intrinsic 
overflow/underflow constraints.
For the particular case of Gauss--Legendre
quadrature  ($\alpha=\beta=0$), a specific asymptotic approach is considered, which yields the most efficient MATLAB implementation available so far. The Gauss--Laguerre and Gauss--Hermite algorithms incorporate subsampling, and scaling is also available in order to extend the computational range. 
Gauss--Radau and Gauss--Lobatto variants are also considered, along with the computation of the associated barycentric weights. Additionally, arbitrary-precision algorithms (in Maple) are offered for the symmetric cases (Gauss--Gegenbauer and Gauss--Hermite), which can be used to compute thousands of nodes with hundreds of digits in a matter of seconds.

\end{abstract}



\section{Introduction}

Gaussian quadrature rules for computing integrals  
\[
I(f)=\int_a^b f(x) w(x)dx\approx \sum_{i=1}^n w_i f(x_i)=Q_n (f),
\]
where $ w(x) $ is a positive weight function defined on the interval $ [a,b] $ (with $a$ and/or $b$ possibly infinite), 
are quadrature rules with maximal polynomial exactness. Specifically, they satisfy $ I(x^k)=Q_n(x^k) $ for $ k=0,1,\dots,2n-1 $. 

As is well known, the nodes $x_i$ of an $n$-point Gaussian quadrature rule correspond to the zeros of the orthogonal polynomial $ p_n(x) $ of degree $ n $ (unique up to a normalization factor) that satisfies $I(x^k p_n(x))=0, \quad k=0,1,\dots,n-1$.  
The quadrature weights are related to the values of the orthogonal polynomials at these nodes. The nodes and weights can also
 be expressed in terms of the eigenvalues and eigenvectors of the Jacobi matrix associated with the three-term recurrence relation satisfied by the family of orthogonal polynomials $\{p_k (x)\}_{k\in \mathbb{N}}$. See, for instance, \cite[sect. 5.3]{Gil:2007:NMSF}.  
Particularly relevant are the so-called classical quadrature rules, for which explicit expressions for the orthogonal polynomials are known. These include Gauss--Jacobi (Gauss--Legendre as a particular case), Gauss--Laguerre and Gauss--Hermite quadratures. 
In these cases, the weights of the $ n $-point rule can be expressed in terms of the derivative of the corresponding orthogonal polynomial evaluated at the nodes. Table 1 summarizes the fundamental relations for these quadratures.  

\begin{table}
\begin{center}
\begin{tabular}{|c|c|c|c|}
\hline
& Hermite & Laguerre & Jacobi\\
&& ($\alpha>-1$) & ($\alpha>-1,\,\beta>-1$)\\
\hline
$(a,b)$ & $(-\infty,+\infty)$         &    $(0,+\infty)$      & $(-1,1)$\\
$w(x)$ & $e^{-x^2}$ & $x^\alpha e^{-x}$& $(1-x)^\alpha (1+x)^\beta$\\
$p_n(x)$ & $H_n(x)$ & $L_n^{(\alpha)}(x)$ & $P_n^{(\alpha,\beta)}(x)$ \\
$k_n$ &    $2^n$   & $(-1)^n/n!$ & $\Frac{(n+\alpha+\beta+1)_n}{2^n n!}$\\
$w_i$ & $\Frac{\sqrt{\pi} 2^{n+1}n!}{H'_{n}(x_i)^2}$ &
$\Frac{\Gamma (n+\alpha+1)}{n! x_i L_n^{(\alpha)\prime} (x_i)^2}$ & 
$\Frac{M_{n,\alpha,\beta}}
{(1-x_i^2) P_n^{(\alpha,\beta)\prime} (x_i)^2}$
\\
\hline
\end{tabular}
\end{center}
\caption{Classical Gaussian quadratures: interval of integration $(a,b)$, weight function $w(x)$, 
 orthogonal polynomial $p_n(x)$ and its leading order coefficient $k_n$ \cite[table 18.3.1]{Koornwinder:2010:OP}, and 
 expression for the weights \cite{Gil:2018:AAT}. The constant $M_{n,\alpha,\beta}$ is $M_{n,\alpha,\beta}=2^{\alpha+\beta+1}\Gamma (n+\alpha+1)\Gamma (n+\beta+1)/(n! \Gamma (n+\alpha+\beta+1))$.
 }
 \label{table1}
\end{table}

Gaussian quadrature rules exhibit fast convergence, particularly for analytic functions \cite{GautschiV}, making them widely used in numerous applications, including spectral methods. This does not imply, however, the superiority of any Gaussian quadrature rule over other quadrature rules. For instance, it has been observed that the simpler Clenshaw-Curtis quadrature rule is competitive with Gauss--Legendre in terms of convergence rate \cite{Trefethen,Tref22,Bornemann}.

On the other hand, the computation of Gauss--Legendre quadratures has dramatically improved, becoming almost elementary in terms of asymptotic approximations \cite{Bogaert:2014:IFC,Opsomer:2023:HOA}, as we later explain. Furthermore, for classical Gauss quadratures, several algorithms have been described that can efficiently compute even millions of nodes \cite{Bogaert:2014:IFC,Bremer,Gil:2018:AAT,Gil:2019:FRA,Gil:2019:NCO,Gil:2021:FAR,Hale:2013:FAC}, clearly overruling the well-known Golub-Welsch algorithm \cite{Golub} (see also \cite[sect. 5.3.2]{Gil:2007:NMSF}).
The algorithms we present in this paper share this capability, while enlarging the range of computation of previous software packages and improving performance. For more general weight functions outside the classical case, these methods are not available, and their computation unavoidably requires slower Golub--Welsch-type approaches (see Gautschi's monograph \cite{GautschiB} and \cite{Lauda} for a more recent example).

Our algorithms for classical Gauss quadratures are based on the asymptotic methods described in \cite{Gil:2018:AAT} 
(Gauss--Hermite and Gauss--Laguerre) and \cite{Gil:2019:NCO} (Gauss--Jacobi), as well as on the globally convergent iterative methods of \cite{Gil:2019:FRA} (Gauss--Hermite and Gauss--Laguerre) and \cite{Gil:2021:FAR} (Gauss--Jacobi). 
The MATLAB and Maple programs implementing these algorithms are available in the GitHub repository: 
\begin{center}
\url{http://github.com/NumericalQuadrature/GaussQuadrature}
\end{center}

In comparing our MATLAB algorithms, we consider the most up-to-date MATLAB implementations developed so far, which, in our experience, are the set of algorithms provided with the Chebfun package. The Julia package FastGaussQuadrature.jl \cite{FastG}, largely influenced by Chebfun, is also a worthy contender.
Adaptations of our MATLAB algorithms to other programming languages will be considered in the future. Regarding the Maple algorithms, there are no current alternatives, as these are, to our knowledge, the first arbitrary-accuracy algorithms for Gauss–Gegenbauer and Gauss–Hermite quadratures.

The asymptotic and iterative approaches used in our algorithms have distinct advantages and their combination, as we demonstrate in this paper, 
yields the fastest and most accurate methods available, with the broadest range of validity with respect to the parameters. For example, in the Gauss--Jacobi case 
our algorithms are not subject to upper restrictions on the parameters $\alpha$ and $\beta$, which both the Chebfun and Julia packages have\footnote{The Chebfun algorithm issues a warning when $\max\{\alpha,\beta\}>5$ while the Julia package uses the slower Golub--Welsch algorithm in that
case.} 
($-1<\alpha,\beta\le 5$); the iterative algorithm, being globally convergent and independent of asymptotic approximations, enables the 
computation for large parameters. 
The iterative method, however, suffers a mild loss of accuracy as the degree increases, and in that case the purely
asymptotic methods are the perfect complement. Similarly, the iterative method allows the computation of Gauss--Laguerre quadratures for almost unbounded
values of the parameter $\alpha$, with asymptotic methods as an ideal replacement for large degrees; 
for Gauss--Laguerre, in addition, scaling can be used to further extend the computational range. Scaling can also be considered in the Gauss--Hermite case. Subsampling  is also supported, both for Gauss--Laguerre and Gauss--Hermite quadratures.

The paper is structured as follows. First, for each of the three classical Gaussian quadratures, we describe the different computational methods (iterative and asymptotic), emphasizing aspects specific to each approach that were not fully detailed in previous publications. We begin with Gauss--Jacobi quadrature, with Gauss--Chebyshev and Gauss--Legendre as particularly important cases. For the Legendre case, a highly efficient method based purely on asymptotics and a look-up table for low degrees is presented. After discussing Gauss--Jacobi quadrature, the Gauss--Laguerre and Gauss--Hermite cases are described. Subsequently, we discuss the computation of Gauss--Radau and Gauss--Lobatto quadratures and the evaluation of barycentric weights for interpolation at orthogonal polynomial nodes. Finally, the last section summarizes the set of MATLAB and Maple algorithms accompanying this work.

\section{Gauss--Jacobi quadrature}

We consider separately the Chebyshev cases $|\alpha|=|\beta|=1/2$, the Legendre case $\alpha=\beta=0$ and 
the more general Jacobi case (including Gauss--Gegenbauer quadrature $\alpha=\beta$). 
Before describing these cases in detail, we
give a brief summary of some properties of the general Gauss--Jacobi case that will be useful in the subsequent discussion.

Jacobi polynomials are solutions of the differential equation
$$
(1-x^2)y''(x)+[\beta-\alpha-(\alpha+\beta+2)x]y'(x)+n(n+\alpha+\beta+1) y(x)=0.
$$
With the change of variable $x=\cos\theta$, the function
\begin{equation}
\label{Yt}
Y(\theta)=\left(1-\cos\theta\right)^{(\alpha +1/2)/2}\left(1+\cos\theta\right)^{(\beta +1/2)/2}P_n^{(\alpha,\beta)}(\cos\theta),
\end{equation}
satisfies the ODE
\begin{equation}
\label{EDOt}
\ddot{Y}(\theta)+\left[\left(n+\frac{\alpha+\beta+1}{2}\right)^2
+\frac{\frac14-\alpha^2}{4\sin^2\frac{\theta}{2}}+\frac{\frac14-\beta^2}{4\cos^2\frac{\theta}{2}}
\right]Y(\theta)=0,
\end{equation}
where dots mean derivative with respect to $\theta$. In this variable, the weights can be written
\begin{equation}
\begin{array}{l}
w_i=M_{n,\alpha,\beta}
(1-x_i)^{\alpha+1/2}(1+x_i)^{\beta+1/2}\dot{Y}(\theta_i)^{-2},
 \end{array}
 \end{equation}
 where $x_i=\cos\theta_i$ are the nodes and the constant $M_{n,\alpha,\beta}$ is the same as in Table \ref{table1}. 
 
 The values $|\dot{Y}(\theta_i)|$ have a small variation as a function of $i$.
 This is a consequence of the fact that the coefficient of the ODE (\ref{EDOt}) is nearly constant as $n\rightarrow +\infty$ in any 
 compact subinterval of $(-1,1)$ (and therefore the solutions are close to trigonometric functions in a 
 first approximation). This is consistent with the known asymptotic approximations as $n\rightarrow +\infty$ 
 (see for instance \cite{Opsomer:2023:HOA}, Eq. (A.8)):
 \begin{equation}
 \label{circle}
 w_i \sim \Frac{\pi}{\kappa} (1-x_i)^{\alpha+1/2}(1+x_i)^{\beta+1/2}(1+{\cal O}(\kappa^{-2})),\kappa\rightarrow +\infty
 \end{equation}
 where 
 \begin{equation}
 \label{kappa}
 \kappa=n+\frac{\alpha+\beta+1}{2}.
 \end{equation}
 and $x_i$ inside any fixed compact subinterval 
 $[a,b]$ of $(-1,1)$.
 
\subsection{Gauss--Chebyshev quadratures}
\label{chebq}

For the Chebyshev cases, $|\alpha|=|\beta|=1/2$, simple explicit formulas exist for the nodes and weights 
\cite[Thm. 8.4]{Mason:2003:CP}. The four cases can be summarized in a single formula as 
follows:
\begin{equation}
\label{chebx}
\begin{array}{l}
x_k=\cos\left((n-k+\frac{\alpha}{2}+\frac34)\Frac{\pi}{\kappa}\right),\, k=1,\ldots n,\\
\\
w_k=\Frac{\pi}{\kappa}(1-x_k)^{\alpha+1/2} (1+x_k)^{\beta+1/2},\, k=1,\ldots n,
\end{array}
\end{equation}
with $\kappa$ as in Eq. (\ref{kappa}). We have chosen the ordering $-1<x_1<...<\ldots x_n<1$. We notice that these expressions 
correspond to the dominant terms in the asymptotic expansions as $\kappa\rightarrow +\infty$ of the general Jacobi nodes and weights (compare, for instance, with Eq. (\ref{circle})); this is no coincidence, as the coefficient of the ODE (\ref{EDOt}) is constant for 
$|\alpha|=|\beta|=1/2$.

Care must be taken in order to avoid loss of accuracy in the evaluation of the trigonometric functions involved in (\ref{chebx}).
 For instance, the previous expression for $x_k$ is not convenient for computing the nodes close to
 zero which, if $\alpha<<n$, correspond to values $k\simeq n/2$. Indeed, we 
 have $x_k=\cos\theta_k$ with $\theta_k= \pi/2+ \delta_k$, with $|\delta_k|$ small; then $x_k=-\sin(\delta_k)\simeq -\delta_k$ and 
 $\delta_k=\theta_k-\pi/2$ suffers loss of accuracy by cancellation. In this case it is preferable to use the expression
 $
 x_k=-\sin((\frac{n}{2}-k+\frac{\alpha-\beta+2}{4})\Frac{\pi}{\kappa})
 $
 where the small number $n/2-k$ (compared to $\kappa$) is computed exactly.
 Similarly precautions must be taken for the factors $1\pm x_k$ appearing in the weights. 
 
 The careful handling of trigonometric functions remains crucial in the general Jacobi case when computing asymptotic expansions in elementary functions.
 
\subsection{Gauss--Legendre quadrature}
\label{gaussle}

The Gauss--Legendre case is a particularly important one and, in fact, it is customary to use the expression ``Gauss quadrature" to
refer to this method alone. Because $x_k=-x_{n-k+1}$ ($-1<x_1<...<x_n<1$) and $w_k=w_{n-k+1}$ we only need to 
compute the non-negative nodes and its corresponding weights.

For Gauss--Legendre quadrature, the Bessel-type expansions of Legendre polynomials \cite{Szego:1934:UEA} enable explicit asymptotic expressions for nodes and weights, as in \cite{Bogaert:2014:IFC}. The Chebfun package {\bf legpts} \cite{Driscoll:2013:chebfun} follows this approach for $n>100$. We also use Bogaert’s expansions, but only for a fixed number of nodes and weights: for $n>80$, the 35 largest positive nodes and corresponding weights $ w_i $ are computed using Bogaert’s formulas, while expansions in elementary functions are considered for the rest. 
For $n\le 80$, we use a look-up table storing precomputed nodes and weights with full double precision ($16$ significant digits), generated using the arbitrary accuracy implementation for Gauss--Gegenbauer quadrature based on the iterative method described in 
\cite{Gil:2021:FAR}, which is provided in the software package accompanying this paper.

In the Bogaert formulas the nodes can be computed as $x_{n-k+1}=\cos\theta_k$, $k=1,\ldots n$, where $\theta_k$ is approximated by
\begin{equation}
\label{bog1}
\theta_k= \displaystyle\sum_{i=0}^{3}F_i (\alpha_k,\cot\alpha_k) \kappa^{-2i}+{\cal O}(\kappa^{-8}),
\end{equation}
and the weights by
\begin{equation}
\label{bog2}
w_{n-k+1}=\Frac{2}{\kappa^2 J_1 (j_k)^2}\Frac{\sin \alpha_k}{\alpha_k}\left(\displaystyle\sum_{i=0}^3 
W_i (\alpha_k,\cot\alpha_k) \kappa^{-2i}+{\cal O}(\kappa^{-8})\right)^{-1},
\end{equation}
where the ${\cal O}(\kappa^{-8})$ terms are not needed in our algorithms. Because only the nonnegative nodes need to be computed, it is enough to take 
$k=1,\ldots \lfloor (n+1)/2\rfloor$.

In these expressions $\alpha_k=j_k/\kappa$, where $j_k$ is
the $k$-th positive zero of the Bessel function of order $0$ ($J_0(x)$), $\kappa =n+1/2$, and $J_1(x)$ is the
Bessel function of order $1$. The functions $F_i$ and $W_i$ 
are given in \cite{Bogaert:2014:IFC} for $i=0,\ldots 5$
and the first two are $F_0 (x,u)=x$, $F_1 (u,x)=(ux-1)/(8x)$,  $W_0(u,x)=1$ and $W_1 (u,x)=(ux+x^2-1)/(8x^2)$. 
For our purpose we only need the first four functions $F_i$ and $G_i$. 

Because we are computing with Bogaert's formulas only for $35$ nodes, we can use precomputed values
for the zeros of the Bessel function $j_k$ and 
the values $J_1 (j_k)$. 
 
With respect to the expansions in terms of elementary functions, the
 expansion described in \cite[sect. 2.2.2]{Gil:2019:NCO},  particularizing for 
 $\alpha=\beta=0$, provides an elementary asymptotic approximation for
 the nodes. An expansion for the weights is obtained by substituting the expansion 
 for the nodes in the
expression for the weights in Table \ref{table1}, considering the expansion for 
$P_n'(x)$ for large $n$ and then re-expanding, similarly
as was done in \cite{Pouso:2024:URB} for the first few terms. 
We skip further details of the computation. The results coincide with
those obtained in \cite{Opsomer:2023:HOA} by using the Riemann-Hilbert approach and
can be written as follows:

\begin{equation}
\label{elem}
\begin{array}{ll}
x_k=& c_k\left(1-\Frac{1}{8\kappa^{2}}+\Frac{33+28r_k^{2}}{384\kappa^{4}}-\Frac{865+2060 r_k^2+1208 r_k^4}{5120\kappa^{6}}
+{\cal O}(\kappa^{-8})\right)\\
& \\
w_k=&\Frac{\pi}{\kappa}s_k\left(1-\Frac{1}{8\kappa^{2}}+\Frac{33+84 r_k^2+56 r_k^4}{384\kappa^{4}}
\right.\\
&\\
& \left.-\Frac{865+6180 r_k^2+10160 r_k^4+4832 r_k^6}{5120\kappa^{6}}+{\cal O}(\kappa^{-8})\right)
\end{array}
\end{equation}

\noindent
where $\kappa=n+1/2$ and, for $k=1,\ldots n$,
\begin{equation}
\label{cos}
\begin{array}{l}
c_k=-\cos\left((4k-1)\Frac{\pi}{4\kappa}\right)=\sin\left((2k-n-1)\Frac{\pi}{2\kappa}\right)\\
s_k=\sin\left((4k-1)\Frac{\pi}{4\kappa}\right),\,r_k=c_k/s_k.
\end{array}
\end{equation}

Observe that $c_k=-c_{n-k+1}$ and $s_k=s_{n-k+1}$, as corresponds to the symmetry relations 
$x_k=-x_{n-k+1}$ and $w_k=w_{n-k+1}$. The exact symmetry relations are not explicit 
in the Bogaert formulas, but hold asymptotically.

As discussed for the Chebyshev case, care has to be taken in the evaluation of the trigonometric functions in order
to avoid accuracy degradation. For instance, the two expressions for $c_k$ in Eq. (\ref{cos}) are not equivalent numerically and the second one is preferable for the nodes close to the origin.

For $n>80$, both expansions (Bessel-type and elementary) coincide within
 full double precision accuracy for the $36^{th}$ largest
node and its corresponding weight; this indicates that 
the elementary expansions suffice to compute all the nodes and weights 
except for the largest $35$ positive nodes. This is confirmed by comparing
 these results with the higher accuracy Maple implementation of Gauss--Gegenbauer quadrature. 
 Additionally, we have compared against the 
Chebfun implementation: our results coincide except that 
for $n<100$ the accuracy of the Chebfun algorithm 
is slightly worse (when asymptotics is not used). Our MATLAB algorithm maintains full double precision for any degree, as the
comparison against the Maple implementation shows.

Summarizing, we consider the following scheme for computing Gauss--Legendre quadratures:

\begin{enumerate}
\item{Compute the non-negative nodes  $x_k$, $k=\lfloor n/2\rfloor+1\ldots n$ and the corresponding weights $w_k$}:
\begin{enumerate}
\item{If} $n\le 80$ use the look-up table.
\item{If} $n>80$ use asymptotic approximations as follows:
\begin{enumerate}
\item{For} the largest $35$ positive nodes and its corresponding weights ($x_{k}$, $w_k$, $k=n-34,\ldots n$) use Bogaert's formulas
(\ref{bog1}) and (\ref{bog2}).
\item{Compute the} rest of non-negative nodes and corresponding weights ($x_k$, $w_k$, $k=\lfloor \frac{n}{2}\rfloor+1\ldots n-35$) with the
elementary expansions (\ref{elem}).
\end{enumerate}
\end{enumerate}
\item{Use} symmetry for the negative nodes and the corresponding weights: $x_k=-x_{n-k+1}$ $w_k=w_{n-k+1}$, $k=1,\ldots \lfloor n/2\rfloor$
\end{enumerate}

\vspace*{0.1cm}

Combining the elementary and Bessel-type expansions has two advantages over using only the latter. First, computing central nodes with the elementary expansion ensures relative, not just absolute, accuracy. This relates to the precise computation of trigonometric functions (see section \ref{chebq}): while straightforward in the elementary case, the Bessel expansion involves $ x_k \sim \cos \theta_{n-k+1} = \cos(j_{n-k+1}/(n+1/2)) $, and if $ \theta_{n-k+1} $ is close to $ \pi/2 $, we must compute $ \theta_{n-k+1} - \pi/2 $ accurately to avoid loss of relative accuracy\footnote{For large enough $n-k$, this is possible via the MacMahon expansion for the zeros.}. Second, the elementary expansion is simpler, avoiding Bessel functions and their zeros, and is highly efficient.

\begin{figure}[htbp]
\label{times}
\begin{center}
 \includegraphics[width=0.8\textwidth]{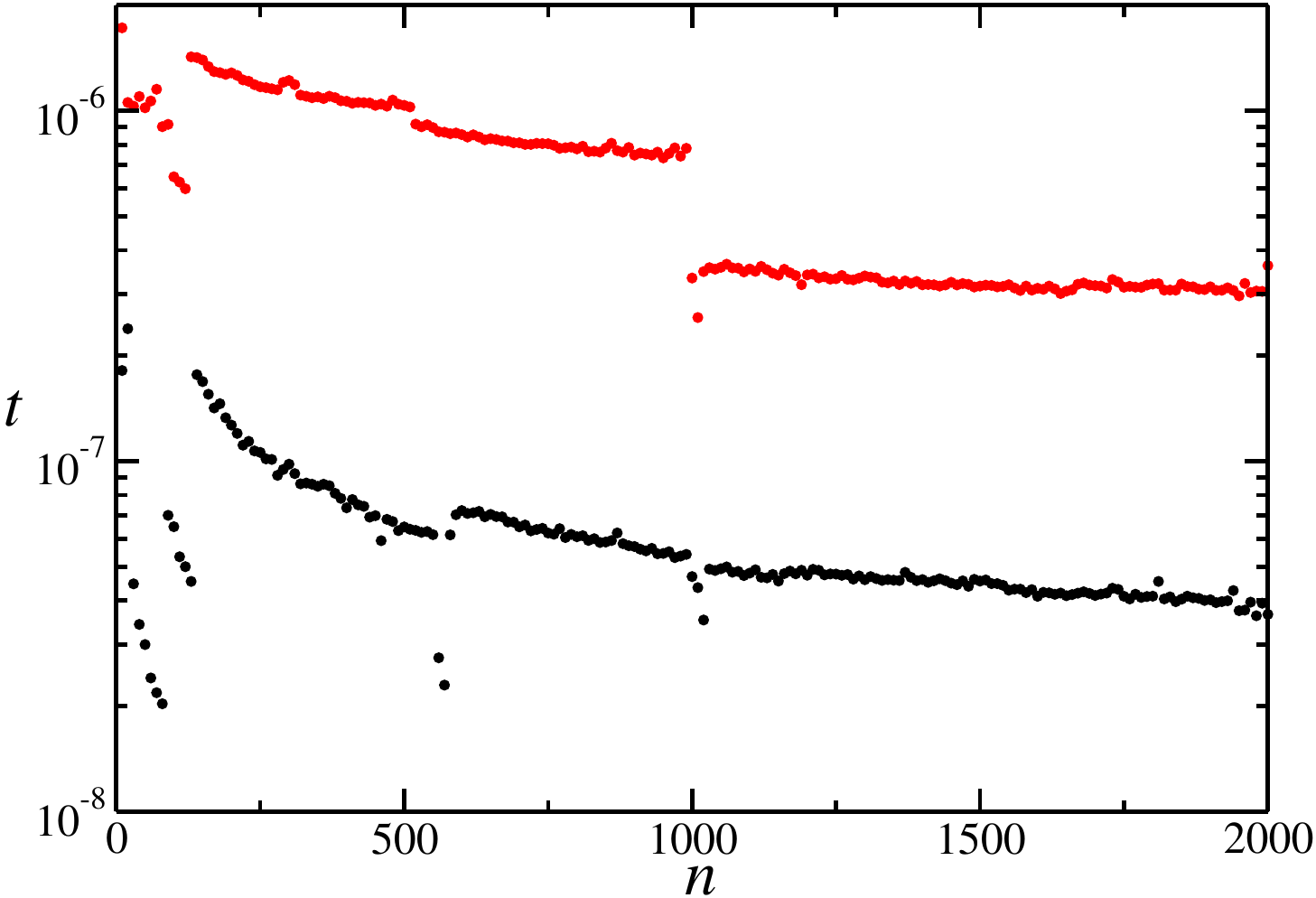} 
\end{center}
 \caption{CPU times in seconds per node and weight spent by the algorithms {\bf legpts} (Chebfun) and {\bf legen} as a function of the 
 degree $n$ (with $n$ from $10$ to $2000$ with step $10$). The times are computed with MATLAB version 2024b running
  on a laptop with  Intel(R) Core(TM) i7-1355U processor. The times are averaged over $100$ evaluations for each value of $n$.}
 \label{cputimes}
\end{figure}

Our algorithm seems to be faster than those based solely on Bessel-type expansions for large enough $n$, as illustrated in Fig. \ref{times}. In this figure, 
the times are measured with the MATLAB instruction {\bf timeit}. Both codes are so simple and efficient that the time spent is too small for a 
reliable measurement, and we have averaged the time over $100$ cycles of computation for each node. Needless to say, the detail of the 
comparative result will depend
on the performance of the computer; the results in Fig. 1 correspond to execution times on a typical laptop, specifically a laptop equipped with and
Intel(R) Core(TM) i7-1355U processor with 10 cores (12 threads).  
The difference can be smaller in higher performance computers. 
For instance, in an Intel(R) Core(TM) i7-14700K with 20 cores (28 threads), we observe roughly a reduction 
by a factor $2$ in the relative difference (not shown), but {\bf legen} is still more that $5$ times faster than {\bf legpts} for $n\le 2000$.

\subsection{General Jacobi case}

The methods considered for the general Gauss--Jacobi case combine the iterative method of \cite{Gil:2021:FAR} with 
the asymptotic methods of \cite{Gil:2019:NCO}. The iterative method provides a fast method
 of computation for
a wide range of parameters; however, the accuracy of the weights worsens gradually as the degree increases 
and the asymptotic methods can be used as a replacement for moderate values of the parameters 
$\alpha$ and $\beta$. 

\subsubsection{Iterative method}

The iterative method used for all the three classical quadratures considered in this paper is based on the fourth order globally convergent method described in \cite{Segura:2010:RCO}. This method applies to homogeneous second order equations in normal
form $\ddot{Y}(z)+\Omega (z)Y(z)=0$ (where dots represent derivatives). 
 Provided the monotonicity properties of $\Omega(z)$ are known, the method computes with certainty and fourth 
order convergence the zeros of any solution
of the ODE, generating monotonically all the zeros in the direction of decreasing $\Omega (z)$.

The main idea is condensed in Theorem 1 of \cite{Gil:2021:FAR}, which we recover here for completeness:

\begin{theorem}
Let $Y(z)$ be a solution of $\ddot{Y}(z)+\Omega (z)Y(z)=0$ and let ${\protect\rm a}$ be 
such that $Y({\protect\rm a} )=0$. 
Let ${\protect\rm b}\neq {\protect\rm a}$ such that $\Omega ({\protect\rm b})>\Omega ({\protect\rm a})>0$ and 
$Y(z)\neq 0$ in the open interval $I$ between ${\protect\rm a}$ and ${\protect\rm b}$. 
Assume that $\Omega(z)$ is differentiable and monotonic in the closure of $I$. Let 
$j=\mbox{sign}({\protect\rm b}-{\protect\rm a})$, then for any $z^{(0)}\in I\bigcup \{{\protect\rm b}\}$, 
the sequence $z^{(i+1)}=T_j(z^{(i)})$, $i=0,1,\ldots$, with
\begin{equation}
T_{j}(z)=z-\Frac{1}{\sqrt{\Omega (z)}}\arctan_{j}\left(\sqrt{\Omega (z)}\Frac{Y(z)}{\dot{Y}(z)}\right)
\end{equation}
and
\begin{equation}
\label{arctan}
\arctan_{j}(\zeta)=\left\{
\begin{array}{l}
\arctan(\zeta)
\mbox{ if }j\zeta > 0,\\ 
\arctan(\zeta)+j\pi \mbox{ if }j\zeta\le 0,\\ 
j\pi/2 \mbox{ if } \zeta=\pm \infty,
\end{array}
\right.
\end{equation}
is such that $\{z^{(i)}\}_{i=1}^{\infty}\subset I$ and it converges monotonically to the root 
${\protect\rm a}$ with order
of convergence $4$ if $\Omega (z)$ is not constant (in which case the method is exact). 
\end{theorem}

The previous theorem gives a procedure
to compute successive zeros in the direction of decreasing values of $\Omega (z)$. 
For a simple application of the method it is convenient that the monotonicity properties of $\Omega (z)$ are
available and that the variation of the coefficient $\Omega (z)$ is as small as possible, because the method is
exact if $\Omega(z)$ is constant. Different changes of variable can be considered such that 
the transformed differential equation becomes suitable in this sense. This has to be combined with an
adequate method for computing the Jacobi polynomial and its derivative (or related functions\footnote{Notice that the fixed point method applies to ODEs in normal form, and that
 transforming the ODE to this form may require a change of function, as for example in Eq. (\ref{Yt})}); 
 in our algorithms, as in
\cite{Gil:2021:FAR}, the main method considered for this purpose is Taylor series,
with the successive derivatives in the series computed by differentiation of the ODE satisfied by Jacobi polynomials 
(or related functions). This approach was earlier considered in \cite{Glaser:2007:AFA}, and it is a simple but effective
method of computation which, however, needs to be complemented with some alternative methods in the Jacobi case (and also
for Gauss--Laguerre, as we later discuss).

 The algorithm for the general Jacobi case combines the use of three different changes of variable. For computing the
 Jacobi polynomial and its derivative (or related functions) with Taylor series the original variable $x$ is used because the differential equation has polynomial
 coefficients and it is therefore easy to differentiate. The principal change of variable considered for the 
 applying the fixed point method is $z=\tanh^{-1}(x)$; the coefficient of the Liouville-transformed ODE in normal form in
 the variable $z$  has even simpler monotonicity properties than that of Eq. (\ref{EDOt}). Finally,
 for nodes close to the end points $\pm 1$ the angular variable $\theta=\arccos x$ is used to refine the results obtained with the fixed
  point method associated to the principal change of variable.
 
 For the particular case of Gauss--Gegenbauer quadrature we have implemented an arbitrary
  precision Maple algorithm in addition to the MATLAB algorithm. The computation of the Gegenbauer polynomial and its derivative 
  (or related
   functions) by Taylor series 
  is particularly simple due to the symmetry of the problem. The angular variable is not used in this Maple algorithm.

 For further details on the iterative method for Gauss--Jacobi quadrature, we refer to \cite{Gil:2021:FAR}.

\subsubsection{Asymptotic method}

Asymptotic methods have been employed for the computation of Gauss--Jacobi quadratures both in combination with iterative methods as in \cite{Hale:2013:FAC} or as independent iteration-free methods
  for Gauss--Legendre \cite{Bogaert:2014:IFC} and Gauss--Jacobi \cite{Gil:2019:NCO}.  
  The Gauss--Hermite and Gauss--Laguerre cases were also described in  \cite{Gil:2018:AAT}. In \cite{Opsomer:2023:HOA}, 
  a Riemann-Hilbert analysis is considered for the three classical cases and some extensions.
 
 In our algorithms, the asymptotic computation of Gauss--Jacobi quadratures is based on the
two types of asymptotic expansions described in \cite{Gil:2019:NCO}: an expansion in terms of elementary functions and
a Bessel-type expansion \cite{Frenzen:1985:AUA}. 
The elementary expansion is
 valid in any fixed compact subinterval of $(-1,1)$ (therefore the
 nodes closest to $\pm 1$ are excluded). Contrarily, the Bessel expansion is uniformly valid in $[-1+\delta,1]$,
 where $\delta$ is an arbitrarily small positive number; because we can use the symmetry property
 $P_n^{(\alpha,\beta)}(x)= (-1)^n P_n^{(\beta,\alpha)}(-x)$ this expansion can be used for computing all the zeros of
 the Jacobi polynomials in $(-1,1)$.

The Bessel-type expansion is the main method of computation in our algorithms. 
 However, the expansion in terms of elementary functions is convenient for computing 
 the nodes with small absolute value  in order to maintain relative accuracy, and not only absolute 
 accuracy, similarly as discussed earlier for Gauss--Legendre quadrature. In our algorithm,
 $10\%$ of the nodes are computed with the 
 elementary expansion and the rest with the Bessel-type expansion.

  We consider enough terms in the asymptotic expansions so that they can be used with nearly double precision accuracy for 
  $n>250$ and moderate $\alpha$ and $\beta$. 
  The range of the parameters $\alpha$ and $\beta$ for which the expansions work accurately 
  is larger as $n$ increases, as we are discussing next.
 
 \subsubsection{Combining the iterative and asymptotic methods}
 
  Our algorithm computes Gauss--Jacobi quadratures using the iterative method and the asymptotic expansions 
  in complementary regions. The iterative method is fast and globally  convergent and can be used
  for practically any values of the parameters, but with some error degradation as the degree increases.
  For that reason, it is convenient to replace the iterative method by asymptotic computation for sufficiently large degrees.
  
  The asymptotic approximations described in \cite{Gil:2019:NCO} hold
  for degrees $n$ much larger that $\alpha$ and $\beta$, and it can be used for larger values of $\alpha$ and 
  $\beta$ as the degree $n$ becomes larger. 
  Asymptotic expansions valid for both large degree and parameters
  were discussed in \cite{Gil:2021:AEO}, but they are not implemented in the current version of the algorithms.

  For deciding which method, iterative or asymptotic, is more accurate,
  we compare them against an alternative method using extended accuracy.  
  We 
  compare against a quadruple precision version of the Golub-Welsch algorithm (GW)  using the Advanpix 
  multiprecision toolbox \cite{advanpix:2023:MCT}.

   We focus on the accuracy of the weights for deciding the regions 
  of application of each method (the computation of the nodes is more accurate).
   We consider three different measures of the error of the weights.
    Denoting by $w_i^{(Q)}$ the quadruple accuracy weights and by $w_i$ the double precision weights, 
   we define the maximum relative error ($\epsilon_{mr}$), average relative error ($\epsilon_{ar}$) and 
   relative total error ($\epsilon_{r}$) as follows:   
   \begin{equation}
   \epsilon_{mr}=\max_{i=1,\ldots n}\left|1-\Frac{w_i}{w_i^{(Q)}}\right|,\,\epsilon_{ar}=\Frac{1}{n}\sum_{i=1}^{n}\left|1-\Frac{w_i}{w_i^{(Q)}}\right|,\,
\epsilon_{rt}=\Frac{\sum_{i=1}^n |w_i-w_i^{(Q)}|}{\sum_{i=1}^n w_i^{(Q)}}.
   \end{equation}
    
    In order to simplify the analysis, we present the accuracy tests in terms of only two 
    parameters: $n$ and $\alpha^2+\beta^2$. Obviously, some information is lost, but 
     this simplification is useful for a simple selection of the methods. 
   
\begin{figure}[htbp]
\begin{minipage}{6.2cm}
\includegraphics[width=\textwidth]{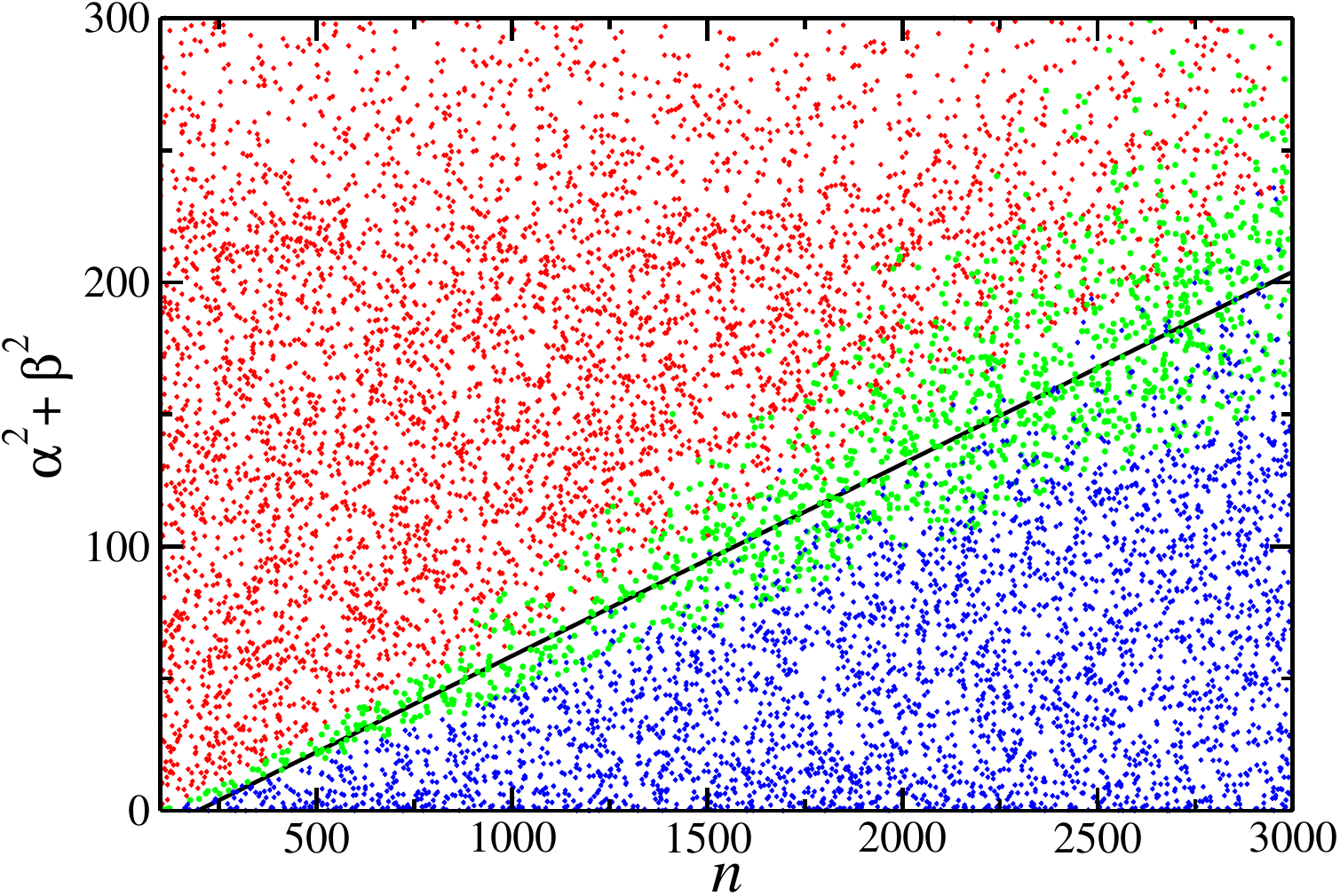}
\end{minipage}
\begin{minipage}{6.2cm}
\includegraphics[width=\textwidth]{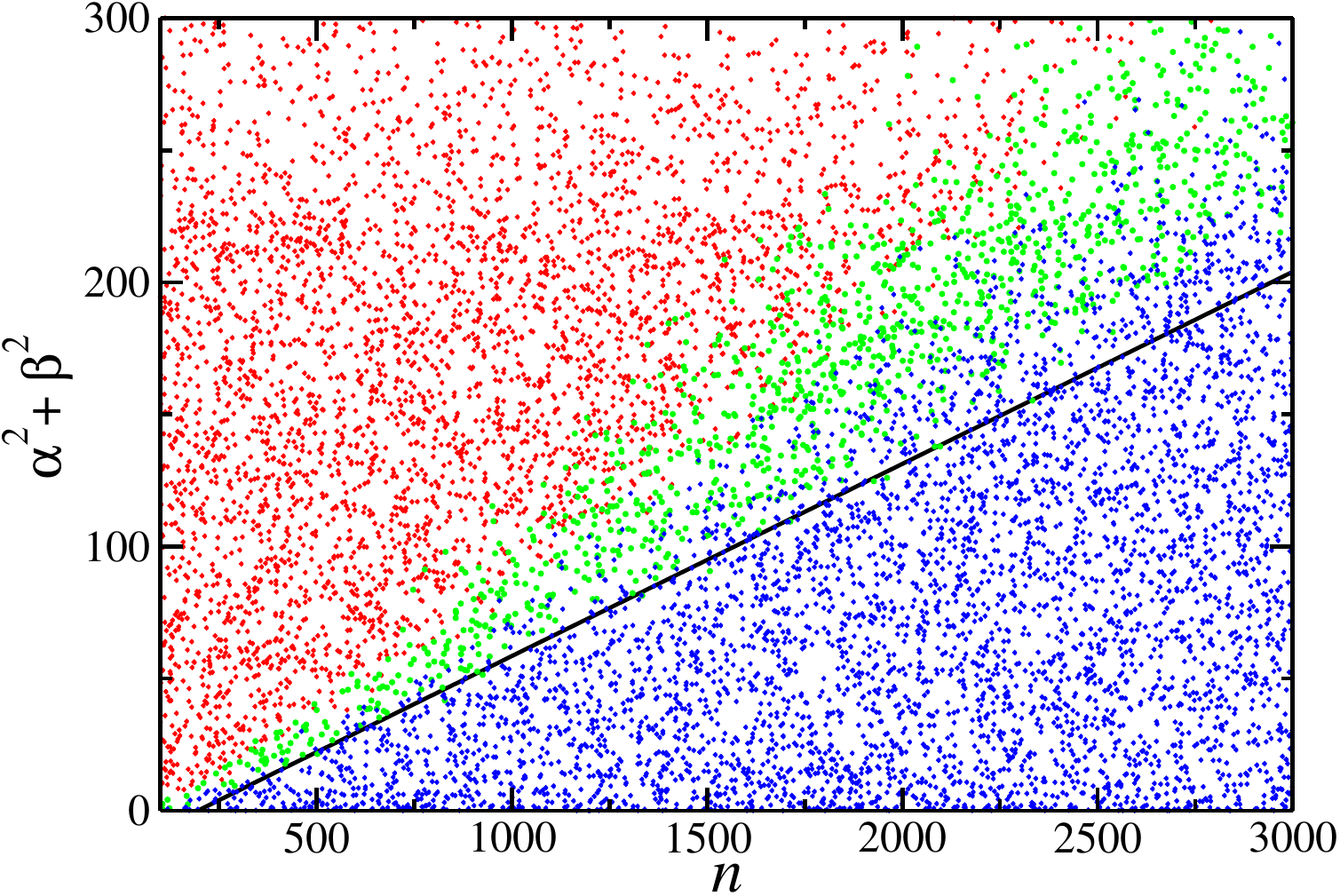}
\end{minipage}
\begin{center}
\includegraphics[width=0.5\textwidth]{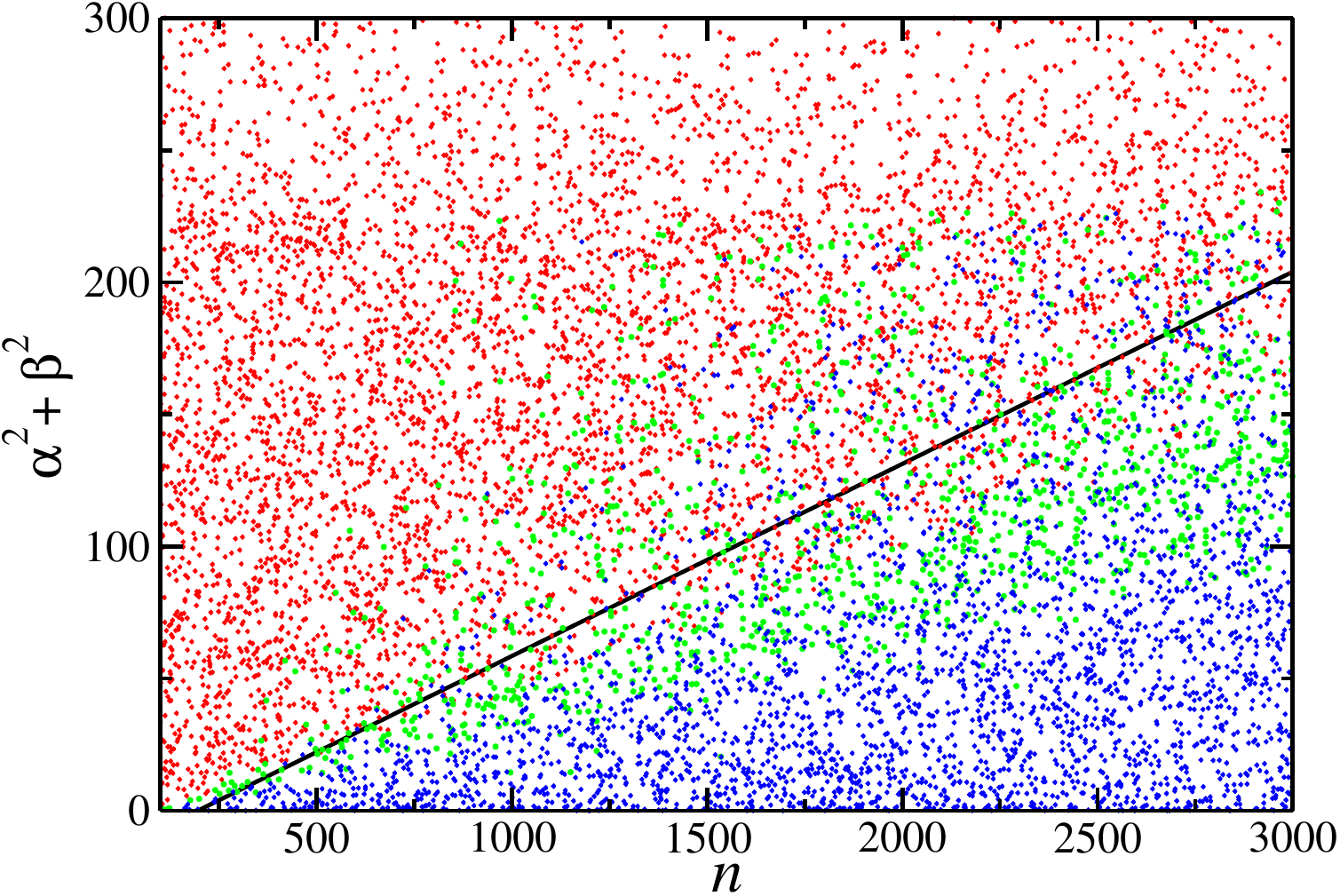}
\end{center}
\caption{
Comparison of the accuracy of the iterative and the asymptotic method for Gauss--Jacobi quadratures, as a function of the degree $n$ and
 $\alpha^2+\beta^2$, 
 for three different measures of accuracy: average
 relative error (left), maximum relative error (center) and relative total error (right). Denoting
 by $R$ 
 the ratio of the accuracy of the asymptotic method over the accuracy of the iterative method, we have $R<1$ at the red points. $1 \le R<10$ at the green points and $R\ge 10$ at the blue points. 
 The separation line 
 $\alpha^2+\beta^2 =\min\{\Frac{4}{55} (n -195)\}$ is also shown. 
}
\label{separa}
\end{figure}

   The accuracy comparison between the asymptotic and iterative methods is shown in Fig. \ref{separa}, where accuracy comparisons are 
   performed for $10^4$ points randomly and uniformly generated in the intervals $n\in [1,3000]$, $\alpha, \beta\in (-1,15]$. 
   The iterative method, as we will see later, turns out to be around a factor $2$ faster than the asymptotic method and, for this reason,
  when the accuracy of both methods is similar we prefer to use the iterative method. 
   With this consideration, and given the accuracy checks shown in Fig. \ref{separa}, a good compromise consist in using the asymptotic 
   method when $n>250$ and $\alpha^2+\beta^2<4(n-195)/55$, and the iterative method otherwise. To this, we
   add the restriction $\alpha^2+\beta^2<2400$ in order to avoid overflow/underflow problems in the computation of the asymptotic expansions (not shown in Fig. \ref{separa}).
 
\subsubsection{Testing the combined algorithm}

We now analyze the performance of the proposed algorithm in terms of computational range, accuracy, and speed. As we demonstrate below, it outperforms previous methods by being faster and by remaining reliable over a wider range of parameter values, including large $\alpha$ and $\beta$.

\begin{paragraph}{Range of computation}

   As shown in Fig. \ref{overf}, our algorithm is able to compute Gauss--Jacobi quadratures parameters for 
   $\alpha,\,\beta \in (-1,1000]$ and $n\le 10^6$, while 
   for larger values of $\alpha$ and $\beta$ the algorithm may fail to produce results.  
   For generating the points for the figure we have considered values of $n$ with 
    $n=\lfloor 10^{6R}\rfloor$  with $R$ uniform pseudorandom  numbers in $[0,1]$; $\alpha$ and $\beta$ are
    also generated by exponentiation of pseudorandom numbers.
   
   In the computable region, some of the weights might be 
   too small when the parameters are large, and those
   weights are set to zero (red points in Fig. \ref{overf}). However, for
   degrees smaller than $1000$ the values of $\alpha$ and $\beta$ may be as large as $100$ and all the nodes and weights
   are computed with some degree of relative accuracy (which we later quantify in detail). 
   For degrees $n\le 10^6$ all the weights given by the algorithm are non-zero for $\alpha$ and $\beta$ smaller than $30$.
   We notice that the region where some of the weights are set as zero (in red color in Fig. \ref{overf}) crosses the 
   separation line between the iterative and the asymptotic regions smoothly, showing that the limitation in the range
   of computable weights is intrinsic of double precision accuracy limitations, and not a particular limitation of the methods.
   
   \begin{figure}[htbp]
\begin{center}
\includegraphics[width=0.6\textwidth]{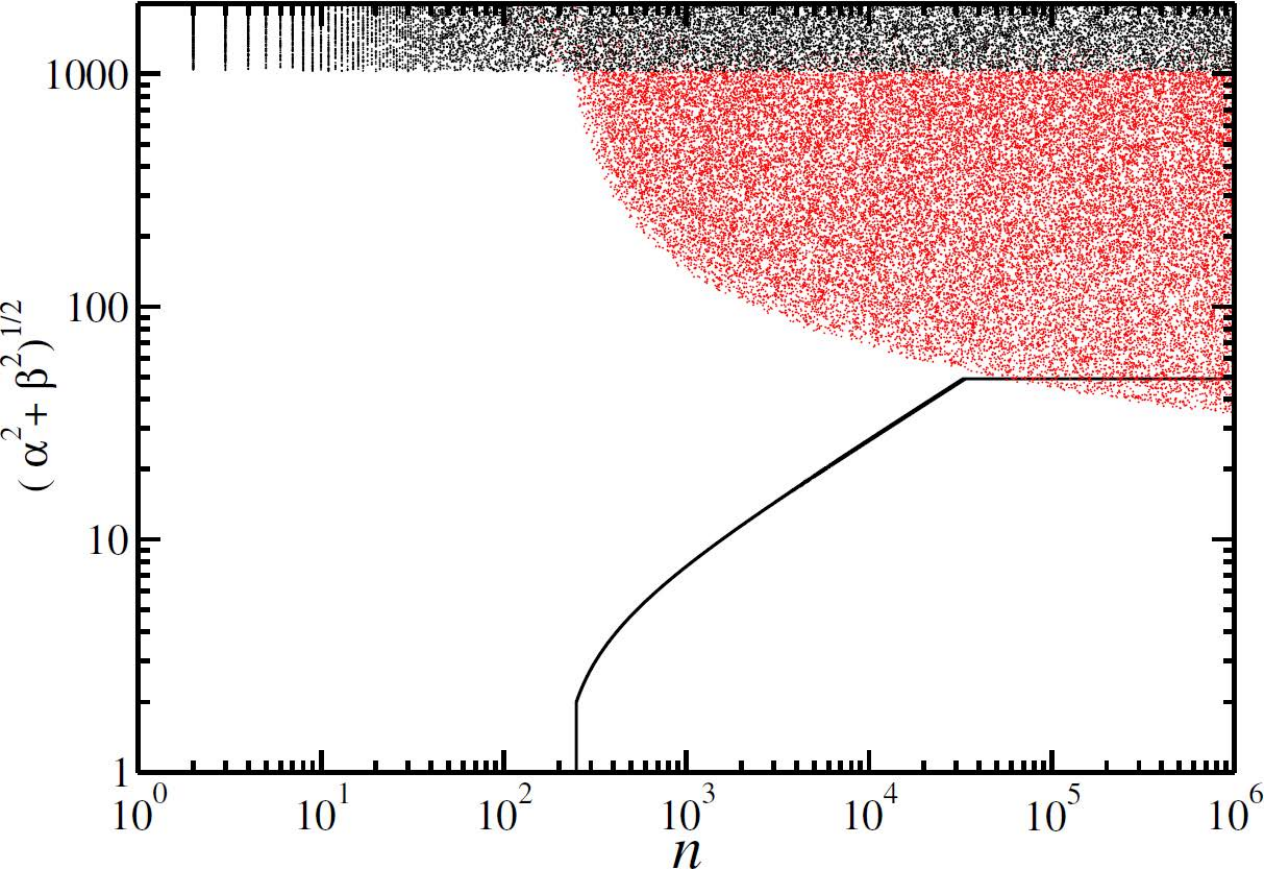}
\end{center}
\caption{Points
where the algorithm fails to compute accurately the nodes and weights.
At the black dots the algorithm is not able to compute any nodes or weights, while at the red nodes some of the weights 
underflow and error degradation is expected for the smallest computed weights. The solid line
  delimits the region where the asymptotic method is applied (bottom-right region in the figure)}
 \label{overf}
\end{figure}

   \end{paragraph}
   
   \begin{paragraph}{Accuracy of the weights}
   In Fig. \ref{precisionw}
   we show the accuracy in the computation of the weights for $n<10^4$, $\alpha,\beta \in (-1,100]$, with $10^4$ points randomly 
   generated by $n=\lfloor 10^{4R}\rfloor$, $\alpha,\, \beta=-2+102^R$, with $R$ denoting different uniform 
   pseudorandom numbers in $[0,1]$. 
We consider two measures of accuracy, the most demanding one (maximum relative accuracy) and the less demanding one (relative total error).
 
The maximum relative error when the iterative method is applied usually takes place at the extremes of the interval,  which corresponds with
 the less significant weights when $\alpha,\beta >-1/2$ (see Eq. (\ref{circle})). 
This is expected because these are the
last computed weights in the iterative algorithm and they are the closest to the singularities of the ODE
(which affects the use of Taylor series). 
For $\alpha,\beta >-1/2$  this error degradation is not of much concern because the extreme weights are the smallest;
however, this is not the case if $\alpha$ and/or $\beta$ are smaller than $-1/2$, when the extreme nodes are the most significant.
For this reason, the algorithm recalculates those extreme nodes and weights in the angular variable $\theta=\arccos x$ for negative values of
$\alpha$ or $\beta$. We refer to \cite[Section 4]{Gil:2021:FAR} for further details.

\begin{figure}[htbp]
\begin{center}
\includegraphics[width=0.49\textwidth]{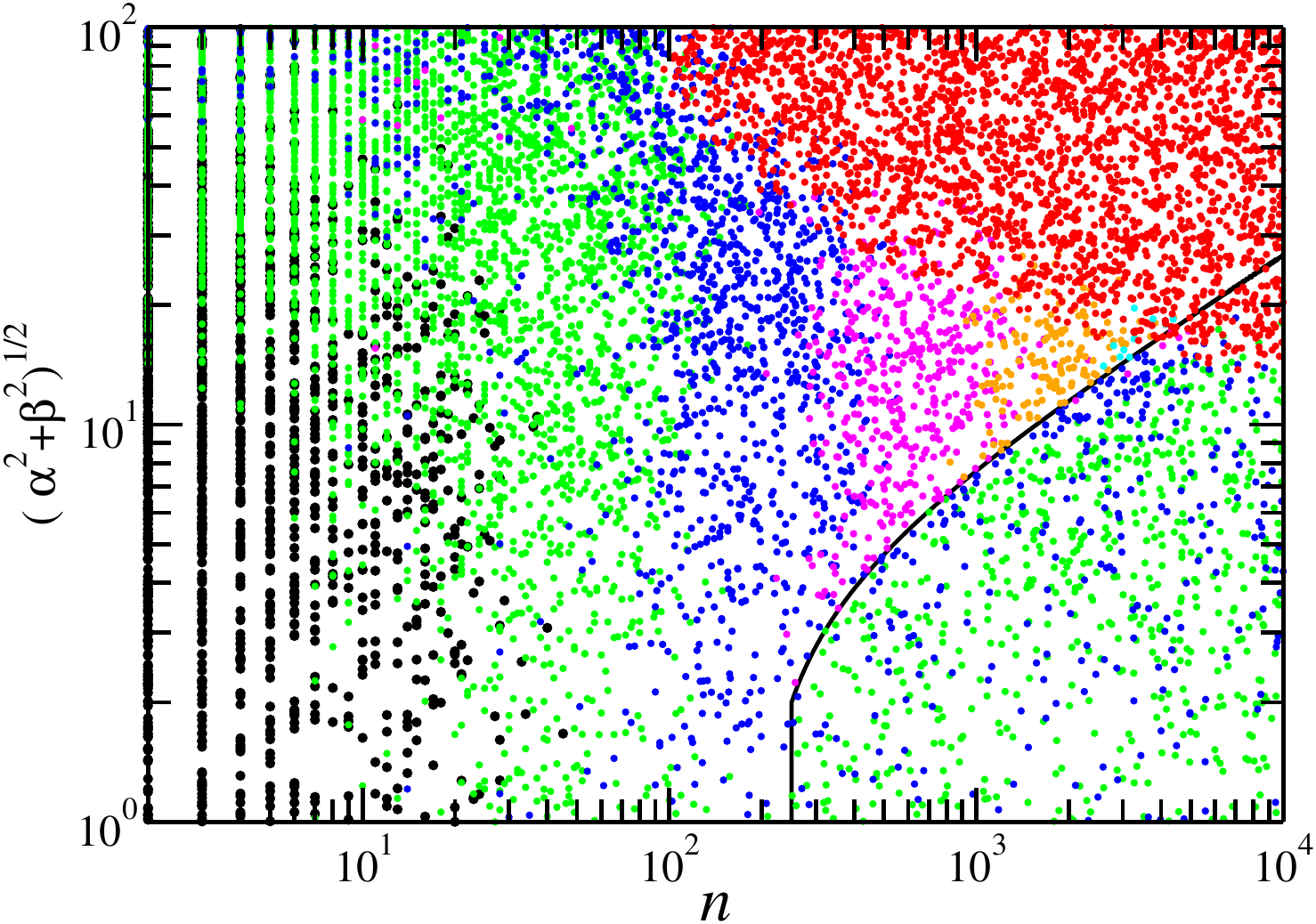}
\includegraphics[width=0.49\textwidth]{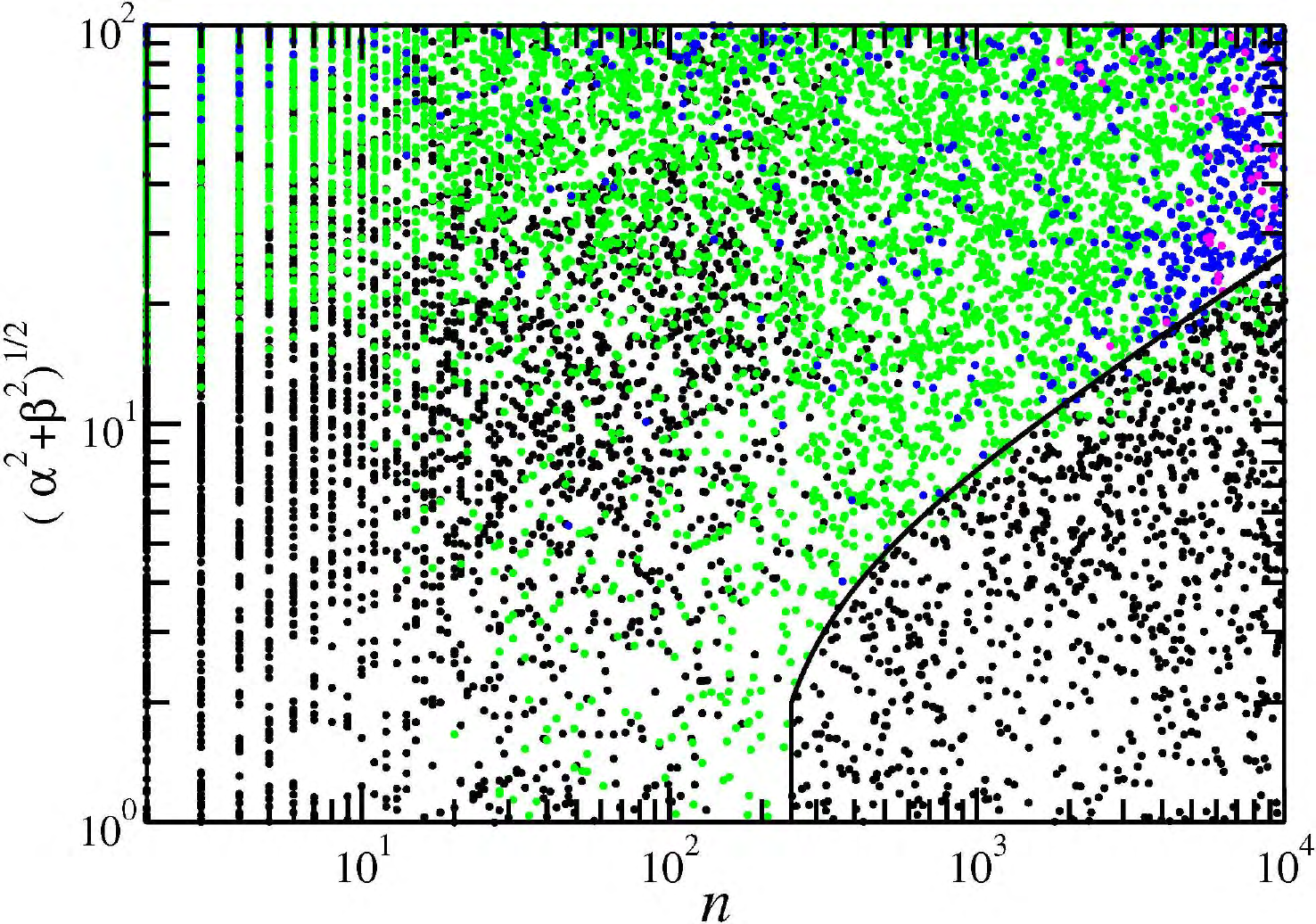}
\end{center}
\caption{Accuracy of the algorithm for computing the weights of Gauss--Jacobi quadrature as functions of $n$ and 
$\sqrt{\alpha^2+\beta^2}$ for two different accuracy measures: maximum relative error (left) and 
relative total error 
(right). The different colors correspond to different error interval:  $(0,10^{-14}]$  for black,
 $(10^{-14},10^{-13}]$ for green, $(10^{-13},10^{-12}]$ for blue, $(10^{-12},10^{-11}]$ for magenta dots,
 $(10^{-11},10^{-10}]$ for orange and 
 $(10^{-10},10^{-9}]$ 
 for cyan. At red dots the GW algorithm fails to compute all the weights. In the right graph, no orange, cyan or red dots
appear.
}
\label{precisionw}
\end{figure}

With respect to the asymptotic computation of the extreme weights it is observed, as
 already commented in \cite{Hale:2013:FAC}, that the computation of Bessel functions 
 with the intrinsic MATLAB programs introduces some loss of precision for some weights corresponding
 to extreme nodes (we observe, as in \cite{Hale:2013:FAC}, that the relative errors are of the order of $10^{-13}$). 
  It should stressed again that this mild loss of relative accuracy is not so relevant as it only affect the smallest
  weights if $\alpha,\beta >-1/2$. 
 \end{paragraph}
  
  \begin{paragraph}{Accuracy of the nodes}
  Our methods produce results with nearly double precision relative accuracy for all the nodes except possibly one. 
   One has to take into account that one node could be exactly zero or close to zero and that relative accuracy
   is not possible for that node. 
   This situation occurs trivially for odd-degree Gauss--Gegenbauer quadrature but also for Jacobi non-symmetrical cases. 
   This is an unavoidable limitation in the computation 
   of the nodes. We notice, however, that we may have this limitation for one node at most, 
   differently to the approach in \cite{Hale:2013:FAC},
   where all the nodes are computed with absolute accuracy. 
   
   In Fig. \ref{precisionx} we test the relative accuracy of the nodes, with points randomly generated in the same way as in Fig.
    \ref{precisionw}. 
   The maximum relative error can be arbitrarily large
    because one of the nodes can be actually zero; in this case
    the second maximum relative error is a more reasonable error measure for the nodes. The second maximum relative error
    is the second greatest value in the set of relative errors for the nodes. We consider both error measures in Fig. \ref{precisionx}.
    
    The plots in Fig. \ref{precisionx}, particularly the right plot, show that the relative accuracy for the nodes is typically smaller
    that $10^{-14}$, with the notable exception of the isolated nodes which are very close to zero. Additionally, we observe
    some loss of accuracy for some values for which the 
    asymptotic method
    is used (blue dots on both plots). This loss of accuracy is due to some lack of accuracy of the MATLAB algorithm for computing Bessel
    functions. 
    The accuracy is in any case always smaller than $10^{-13}$. 
      
\begin{figure}[htbp]
\begin{center}
\includegraphics[width=0.49\textwidth]{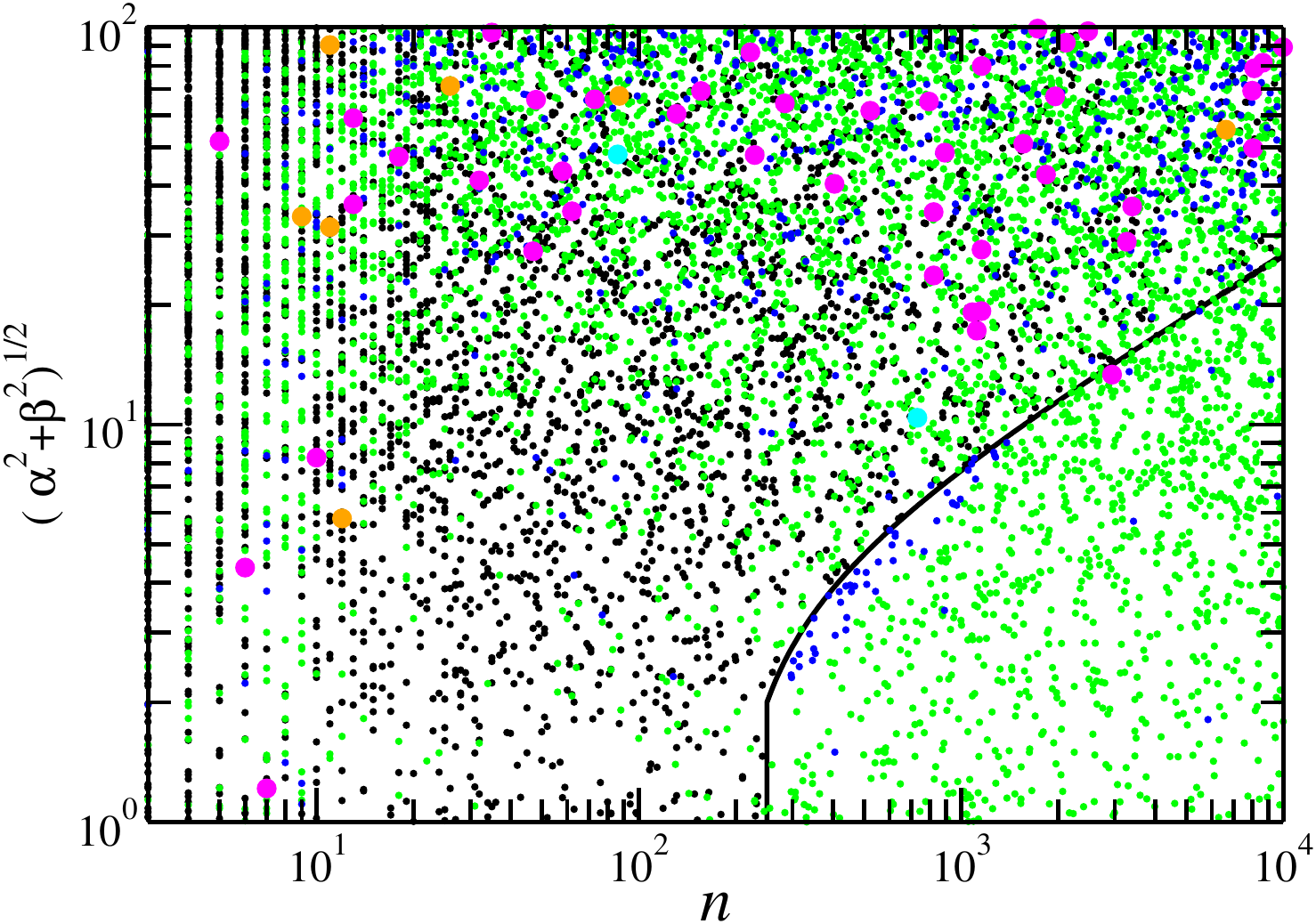}
\includegraphics[width=0.49\textwidth]{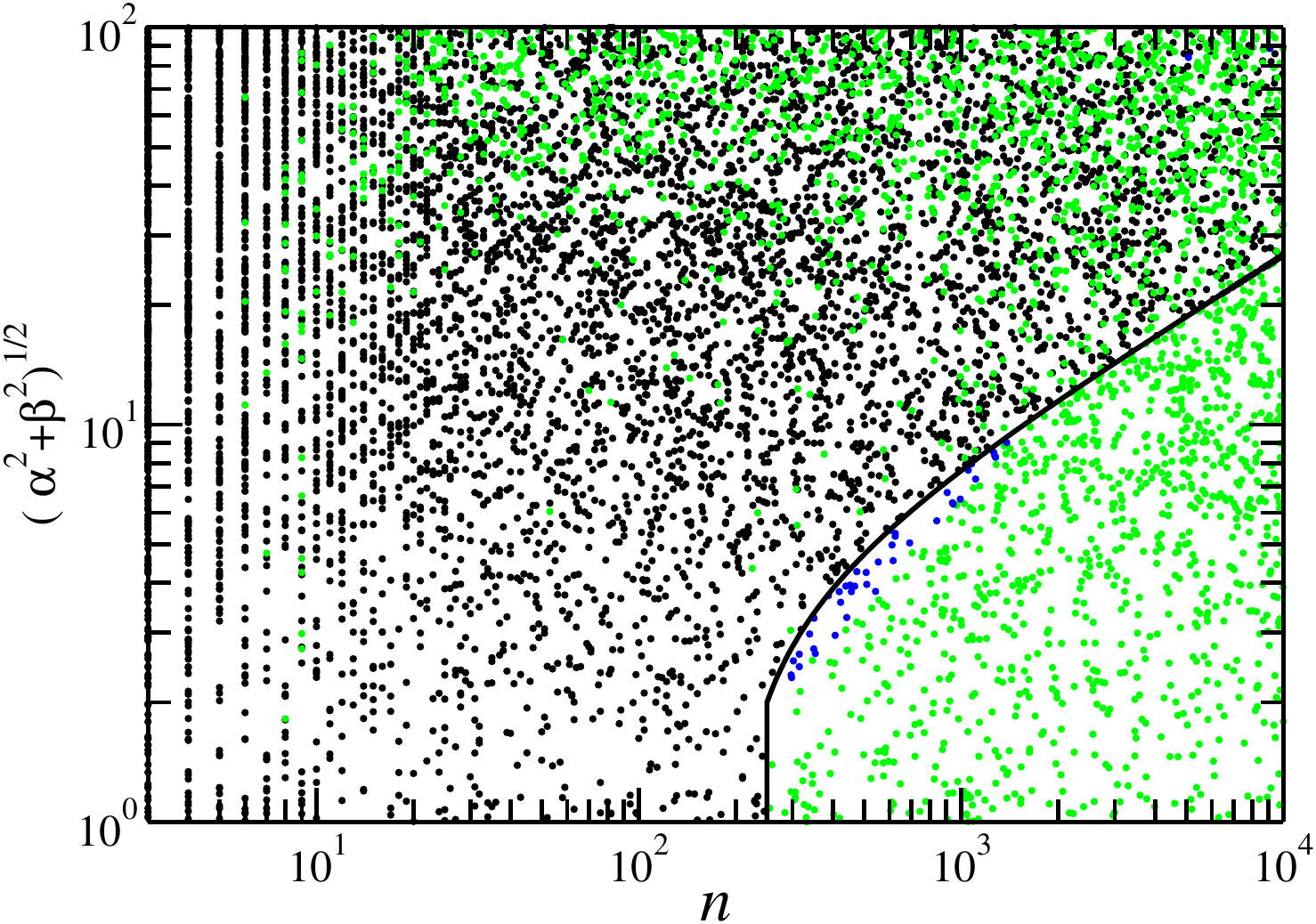}
\end{center}
\caption{Accuracy of the algorithm for computing the nodes of Gauss--Jacobi quadrature 
for two differrent accuracy measures: maximum relative error (left) and 
second maximum relative error
(right). The different colors correspond to different error intervals: $(0,10^{-15}]$  for black,
 $(10^{-15},10^{-14}]$ for green, $(10^{-14},10^{-13}]$ for blue, $[10^{-13},10^{-12})$ for magenta,
 $[10^{-12},10^{-11})$ for orange and 
 $[10^{-11},10^{-10})$ 
 for cyan. Magenta, orange and cyan dots are represented with wider circles for visibility. 
}
\label{precisionx}
\end{figure}

\end{paragraph}   
   
   \begin{paragraph}{Speed of the methods}
   The iterative method is faster than the asymptotic method by around a factor two, and this is why 
   we use the iterative method when both give similar accuracy. On the other hand, our asymptotic method is faster than methods described
   in \cite{Hale:2013:FAC,Driscoll:2013:chebfun} because our expansions (particularly for the nodes) are sufficiently accurate so 
   the iterative refinement is not needed. As a result, our algorithm appears to be significantly faster than previous approaches, as 
   Fig. \ref{cputimes2} shows.

\begin{figure}[htbp]
\begin{center}
\includegraphics[width=0.49\textwidth]{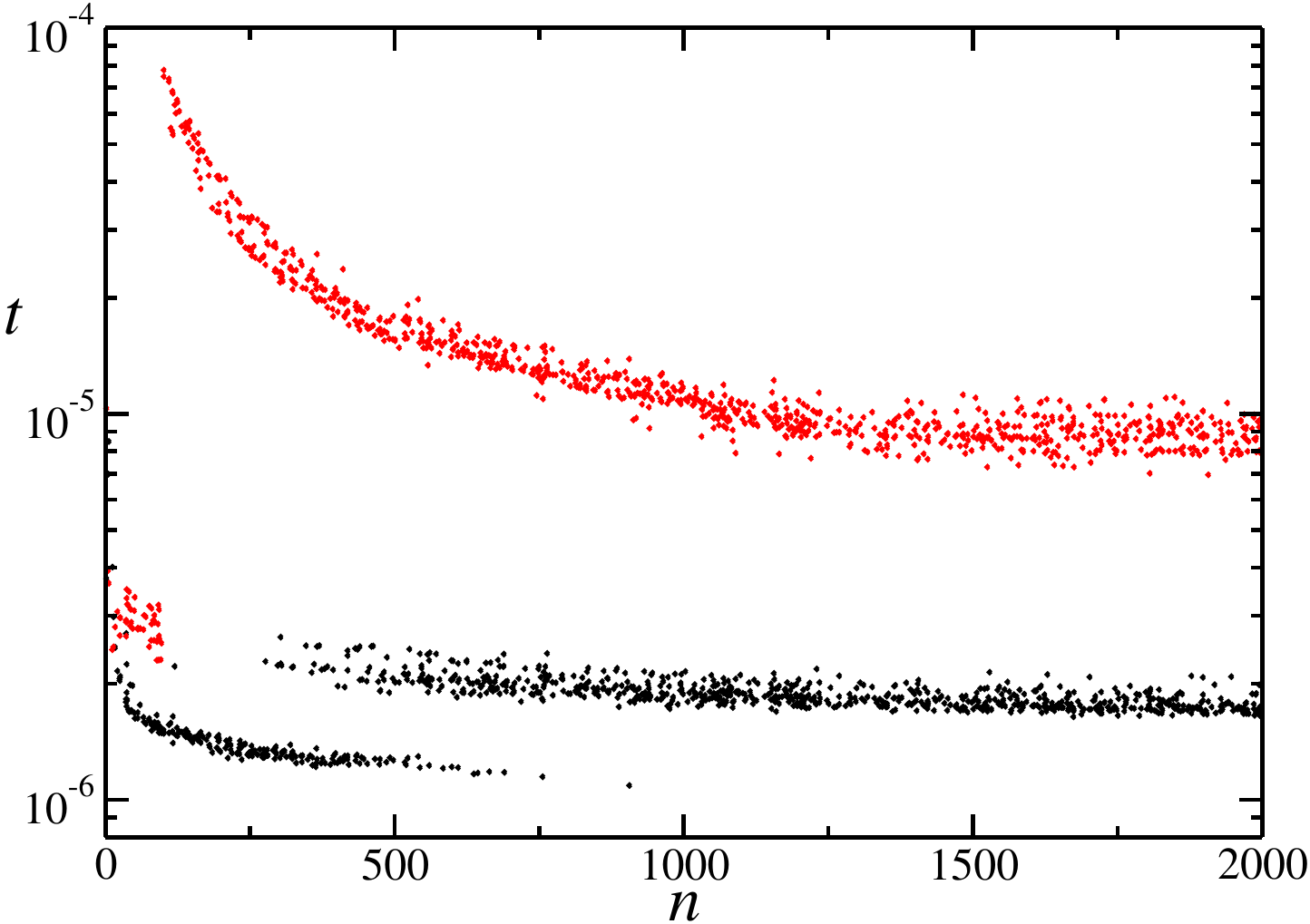}
\includegraphics[width=0.49\textwidth]{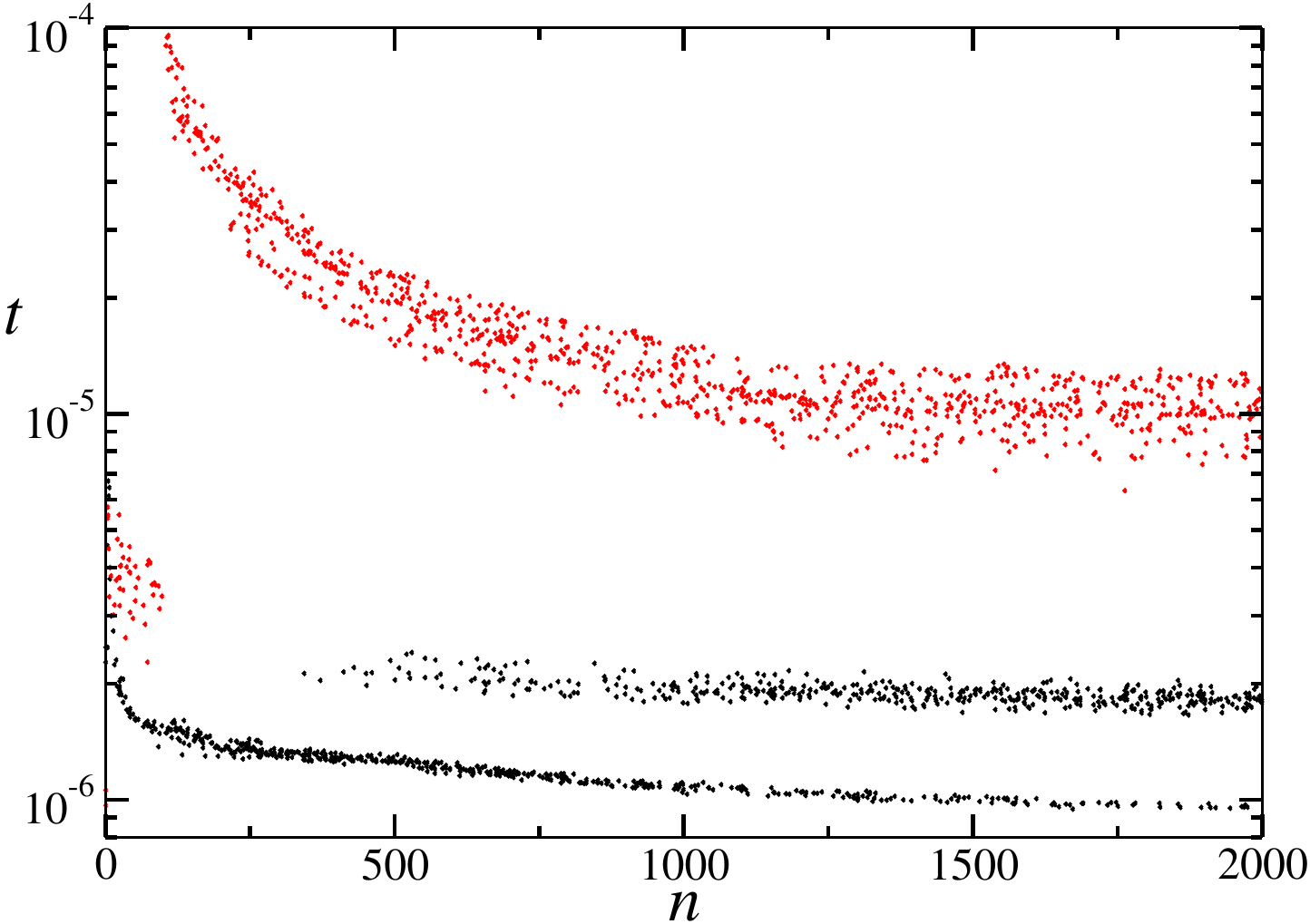}
\end{center}
\caption{CPU times in seconds per node and weight spent by the algorithms {\bf jacpts} (red dots) and {\bf GJ} (black dots) as a function of the 
 degree $n$. The times are computed with MATLAB version 2024b running
  in a laptop with  Intel(R) Core(TM) i7-1355U processor.}
\label{cputimes2}
\end{figure}   
    
    The plots in Fig. \ref{cputimes2} correspond to $10^4$ sets of uniformly generated 
    values of $n$, $\alpha$, and $\beta$, with values of $n$ in the interval $[1,2000]$ and 
    values of $\alpha$ and $\beta$ in $(-1,5]$ in the left figure and in $(-1,10]$ 
    in the right figure. The CPU times are represented as a function of $n$ and, therefore, for each $n$ the points may
    corresponding to different values of $\alpha$ and $\beta$ and then to different methods of computation. That is the case 
    of the black dots, corresponding to our algorithm: the lower curve of black dots are points for which the iterative method
    is applied while the upper curve corresponds to the asymptotic method. We notice that, as advanced, the iterative method is approximately 
    two times faster than the asymptotic method. In the plot on the right, because the range for $\alpha$ and $\beta$ is larger, there are 
    more points where the iterative method is applied.
    
    The reason why the iterative method is so fast lies on two facts: first, the iterative method has fourth order convergence and, second,
    the computations involve relatively simple expressions, in contrast with the more involved expressions for the coefficients in the
    asymptotic series; in addition, one of the asymptotic expansions requires the computation of Bessel function.
    As we discussed, the main reason for using asymptotics is accuracy, particularly for large degrees.
 \end{paragraph}

\section{Gauss--Laguerre quadrature}

As for the Jacobi case, the algorithms for Gauss--Laguerre quadratures combine purely iterative and purely
asymptotic methods, applying each of them in complementary regions. The iterative method is described in
\cite{Gil:2019:FRA} and the asymptotic approximations are developed in \cite{Gil:2018:AAT}. In this section,
 we outline the main ingredients of the iterative method 
 and
 we describe further details of the asymptotic method that were not fully developed in \cite{Gil:2018:AAT}
 and are necessary in the construction of the algorithms. 
 Before this we recall some basic properties.
 
 Laguerre polynomials are solutions of the differential equation
 $$
 x w''(x)+(\alpha+1-x)w'(x)+n w(x)=0.
 $$
 Considering the change $x=z^2$ we have that
 \begin{equation}
\label{yzl}
y(z)=z^{\alpha+1/2} e^{-z^2/2} L_n^{(\alpha)}(z^2)
\end{equation}
satisfies an equation in normal form
\begin{equation}
\label{normL}
\begin{array}{c}
 \ddot{y}(z)+A(z)y(z)=0,\\
 \\
 A(z(x))=-x+\nu+\frac{\frac14- \alpha^2}{x},\,\nu=4n+2\alpha+2.
 \end{array}
\end{equation}

 The Gauss--Laguerre weights (see Table \ref{table1}) can be written
in terms of (\ref{yzl}) as
\begin{equation}
\label{wasym2}
w_i = x_i^{\alpha+\frac12}e^{-x_i} \omega_i ,\,\omega_i=\Frac{4\Gamma (n+\alpha+1)}{n! \left[\dot{y}(z_i)\right]^2},
\end{equation}
where $z_i=\sqrt{x_i}$.
We call the quantities $\omega_i$ scaled weights. The scaled weights have slow variation as a function
of $i$ compared with the unscaled weights $w_i$. This is in accordance with the asymptotic 
estimate for the weights \cite[15.3.15]{Szego:1975:OP}
\begin{equation}
\label{asywl}
w_i\sim \Frac{\pi}{\sqrt{n}}x_i^{\alpha +\frac12} e^{-x_i},\,n\rightarrow \infty
\end{equation}
for the nodes in any fixed closed subinterval of ${\mathbb R}^+$. Comparing this with (\ref{wasym2}) we conclude that $|\dot{y}(z_i)|$
has small variation as a function of $i$ for large $n$, as could be expected because for large $n$ 
the coefficient $A(z)$ of (\ref{normL}) has slow variation in most of the
interval where $A(z)>0$ (where almost all the zeros lie, if not all of them).

\subsection{Iterative method}

For the Gauss--Laguerre case, the iterative method is the main procedure for computing 
the nodes and weights. We again employ the method for second order ODEs described in \cite{Segura:2010:RCO}, using the ODE (\ref{normL}). All the computations necessary can be performed in the  variable
$z$ corresponding to the ODE (\ref{normL}), including local Taylor series for computing $y(z)$ and its derivative, and no
other variables are required (differently from the Gauss--Jacobi case).

The ODE (\ref{normL}) is optimal for several reasons. In the first place, the coefficient $A(z)$ has
simple monotonicity properties, which is important because the fixed point method sweep of \cite{Segura:2010:RCO}
should be applied in the direction of decreasing $A(z)$.  On the other hand, for large $n$
the coefficient $A(z)$ has small variation in most of the oscillation region, and because the nonlinear
method of  \cite{Segura:2010:RCO} is exact for ODEs with constant coefficients, the speed of 
convergence for computing each node tends to improve as $n$ increases. In the third place, the scaled weights $\omega_i$ are
well conditioned as a function of $z_i$ because $\ddot{y}(z_i)=0$. Finally, the resulting method 
admits subsampling in a natural way and one can easily compute only the weights that are greater than a given threshold.

With respect to subsampling, we notice that the main dependence of the weights is carried by the
factor $f(x)=x^{\alpha + 1/2}e^{-x}$ (see (\ref{asywl})), and that $f(x)$ 
has its maximum at $x_M=\alpha+1/2$ when
$\alpha>-1/2$, which is a value close to the starting point for the fixed point method when $|\alpha|>1/2$ ($x_e=\sqrt{\alpha^2-1/4}$). Therefore, for $|\alpha |>1/2$ the fixed point
method computes nodes in the direction of decreasing weights, starting from the most significant
weights (close to $x_e$) and computing the nodes at the right of $x_e$ in increasing order and at the left 
in decreasing order (in both cases in decreasing order for the weights).
For $|\alpha|<1/2$ this is also true because the first computed node is the smallest, which 
corresponds to the largest weight.

 Let us then assume that the 
 first computed node is the $s$-th positive node. After this the algorithm proceeds in the direction of
 decreasing weights, and the ratio $F_i=f(x_i)/f(x_s)=\exp(F_i)$, where
 \begin{equation}
 \label{efe}
 F_i=x_s-x_i+(\alpha+1/2)\log(x_i/x_s)
 \end{equation}
  controls the size of the generated weights $w_i$ with respect
 to the dominant weight $w_s$. If we want to compute only those weights such that $w_i/w_s>\delta$ we can
 stop the algorithm when the index $i$ is such that $F_i<\log(\delta)$. Typically, there are two such values
 of $i$ for $\alpha$ large, one smaller and one large than $s$. Subsampling can save a considerable amount
 of time, particularly for large degrees, because the weights decay exponentially, as (\ref{asywl}) shows.
 
 Because the sum of Laguerre weights is $\sum_{i=1}^n w_i =\Gamma (\alpha+1)$ for large $\alpha$ some weights
 may overflow. For this reason, our algorithms give the option of computing weights with unitary normalization 
 (so that their sum is $1$). The relation between scaled and unscaled weights is the same independently of the
 normalization of the weights, namely:
 $$
 \omega_i= \exp(F_i) w_i.
 $$
 Scaling is considered with respect to the weight $w_s$ described before, in the sense that $w_s=\omega_s$. The
 algorithm also gives as an output the value of the node $x_s$ in order to connect scaled and unscaled weights.
 
 The algorithm proceeds as follows for computing the weights: first the scaled weights with an arbitrary normalization are computed
 and, from those, unscaled weights with this same normalization are obtained. During the computation, the algorithm checks 
 the condition $F_i<\log(\delta)$ if subsampling is required. After this, the weights are normalized to $1$ and 
 this same normalization is applied to the scaled weights. Finally, if the unitary normalization is not chosen, the
 scaled and unscaled weights are multiplied by $\Gamma (\alpha+1)$.

 The scaled weights, particularly with the unitary normalization, can be computed without practical restrictions
 on the degree $n$ and the parameter $\alpha>-1$, with only some mild loss of accuracy as the degree increases. 
 That the algorithm is practical unrestricted is a distinctive characteristic not shared by any previous algorithm.

\subsection{Asymptotic methods}

Our asymptotic algorithm is based on the expansions 
described in \cite{Gil:2018:AAT}. Bessel-type expansions are considered for $80\%$ of the nodes and weights, while the remaining 
$20\%$ (corresponding to the largest nodes) 
are computed with Airy-type expansions. We describe additional details that were
not given in \cite{Gil:2018:AAT} and which are important in the construction of the algorithms.

\subsubsection{Airy-type expansion}

The asymptotic 
computation of the nodes with the Airy expansion is described in detail in \cite[Section 3.2.1]{Gil:2019:NCO}. The starting point is the following representation for $L_n^{(\alpha)}(\nu x)$ discussed in section 3.2 of \cite{Gil:2018:AAT}
(see also \cite{Frenzen:1988:UAE})
\begin{equation}\label{eq:lagder01}
L_n^{(\alpha)}(\nu \tilde{x})=(-1)^n
\frac{ e^{\frac12\nu \tilde{x}}\chi(\zeta)}{2^\alpha\nu^{\frac13}}\left(\Ai\left(\nu^{2/3} \zeta\right)A(\zeta)
+\nu^{-\frac43}\Ai^{\prime}\left(\nu^{2/3}\zeta\right)
B(\zeta)\right)
\end{equation}
with  expansions
$A(\zeta)\sim\sum_{j=0}^\infty\frac{\alpha_{2j}}{\nu^{2j}},\quad B(\zeta)\sim\sum_{j=0}^\infty\frac{\beta_{2j+1}}{\nu^{2j}},\quad \nu \to\infty$,
uniformly for bounded $\alpha$ and  $\tilde{x}\in(\tilde{x}_0,\infty]$, where 
$\tilde{x}_0\in(0,1)$, a fixed number. Details on the coefficients are given in \cite{Gil:2018:AAT}. 
The factor $\chi (\zeta)$ is $\chi(\zeta)=2^{\frac12}\tilde{x}^{-\frac14-\frac12\alpha}\left(\zeta/(\tilde{x}-1)\right)^{\frac14}$.

In the oscillatory
region, where the nodes lie, we have $0<\tilde{x}<1$ and
\begin{equation} 
\label{cambioa}
\dsp{\tfrac23(-\zeta)^{\frac32}=\tfrac12\left(\arccos\sqrt{{\tilde{x}}}-\sqrt{{\tilde{x}-\tilde{x}^2}}\right)}.
\end{equation}

The asymptotic approximation of the nodes is obtained  by inverting the asymptotic series (\ref{eq:lagder01}), the 
first order approximations being the values $\zeta_i^{(0)}$ such that $\Ai\left(\nu^{2/3} \zeta_i^{(0)}\right)=0$, that is $\zeta_i^{(0)}=\nu^{-2/3}a_i$, with 
$a_i$ the $i$-the real (negative) zero of $\Ai (z)$, and the asymptotic series for the nodes in the variable $\zeta$ is written 
\begin{equation}
\label{asyze}
\zeta_i =
\zeta_i^{(0)}+\Frac{\zeta_i^{(1)}}{\nu^2}
+\Frac{\zeta_i^{(2)}}{\nu^4}+\cdots.
\end{equation} 
Details on the computation of the coefficients of this asymptotic expansion are given in 
 \cite{Gil:2018:AAT}. We only mention that for computing the expansion we need the relation (\ref{cambioa}) and its inverse. In our algorithms
 we consider the change $\tilde{x}=\cos^2\Frac{\theta}{2}$, which leads to $\frac83(-\zeta)^{3/2}=\theta-\sin\theta$, which is solved numerically with 
 the Newton method
 \footnote{This change is also considered in \cite{Gil:2018:AAT}, after Eq. (75),
 but it is erroneously written $x=\cos^2 \theta$}. For small $(-\zeta)^{3/2}$ (which is the case for the largest nodes) 
 loss of accuracy is expected due to cancellations, 
 and it is better to write the solution of of $p^3 /6=\theta-\sin\theta$ as 
 $$
 \theta=p\left(1+\Frac{1}{60}p^2+\Frac{1}{1400}p^4+\Frac{1}{25200}p^6+\cdots \right),
 $$
 which, as we have checked, gives full double precision accuracy with $8$ terms if $p<0.45$.

The asymptotic computation of the weights involves 
the function $\dot{y}(z_i)$ in 
(\ref{wasym2}), which requires the derivative of the Laguerre polynomial. In \cite{Gil:2018:AAT}, it is suggested to use the 
property $L_{n}^{(\alpha)\prime}(x)=L_n^{(\alpha)}(x)-L_n^{(\alpha+1)}(x)$ (an alternative is $L_{n}^{(\alpha)\prime}(x)=L_{n-1}^{(\alpha+1)}(x)$).
However, the use of such relations involves the application of (\ref{eq:lagder01}) with an
 argument of the Airy functions different from that which appeared in the computation of the nodes (because
$n$ and/or $\alpha$ are different). We have found that the computation of the Airy functions involved 
is simpler and more accurate by considering an expansion for the derivative with 
exactly the same dependencies.
For this reason, next we obtain the Airy-type expansion for the derivative of the 
Laguerre polynomial.

Differentiating \eqref{eq:lagder01} with respect to $\zeta$ we get
\begin{equation}\label{eq:lagder08}
\begin{array}{ll}
\dsp{\frac{d}{d\zeta}L_n^{(\alpha)}(\nu \tilde{x})=\nu \frac{d\tilde{x}}{d\zeta}{L_n^{(\alpha)}}^\prime(\nu \tilde{x})
=\left(\tfrac12\nu\frac{d\tilde{x}}{d\zeta}+\frac{\chi^\prime(\zeta)}{\chi(\zeta)}\right)L_n^{(\alpha)}(\nu \tilde{x})\ +}\\[8pt]
\quad\quad
\dsp{(-1)^n
\frac{ e^{\frac12\nu \tilde{x}}\chi(\zeta)}{2^\alpha\nu^{\frac13}}
\Biggl(\nu^{2/3}\Ai^\prime\left(\nu^{2/3} \zeta\right)A(\zeta)+\Ai\left(\nu^{2/3} \zeta\right)A^\prime(\zeta)\ +}\\[8pt]
\quad\quad
\dsp{\zeta\Ai\left(\nu^{2/3}\zeta\right)B(\zeta)+\nu^{-4/3}\Ai^\prime\left(\nu^{2/3}\zeta\right)B^\prime(\zeta)\Biggr)},
\end{array}
\end{equation}
which gives
\begin{equation}\label{eq:lagder10}
{L_n^{(\alpha)}}^\prime(\nu \tilde{x})=(-1)^n\frac{ e^{\frac12\nu \tilde{x}}\chi(\zeta)}{2^\alpha\nu^{2/3}}\frac{d\zeta}{d\tilde{x}}
\left(\nu^{-\frac23}\Ai\left(\nu^{2/3} \zeta\right)C(\zeta)
+\Ai^{\prime}\left(\nu^{2/3}\zeta\right)D(\zeta)\right),
\end{equation}
with
\begin{equation}\label{eq:lagder11}
\begin{array}{ll}
\dsp{C(\zeta)= \phi(\zeta)A(\zeta) +A^\prime(\zeta)+\zeta B(\zeta), }\,\\[8pt]
\dsp{D(\zeta)= A(\zeta)+B^\prime(\zeta)\nu^{-2}+\phi(\zeta)B(\zeta)\nu^{-2},  }
\\[8pt]
\dsp{\phi(\zeta)= \tfrac12\nu\frac{d\tilde{x}}{d\zeta}+\frac{\chi^\prime(\zeta)}{\chi(\zeta)}.  }
\end{array}
\end{equation}

We have the expansions
\begin{equation}\label{eq:lagder12}
C(\zeta)\sim\sum_{j=0}^\infty\frac{\gamma_{2j}}{\nu^{2j}},\quad D(\zeta)\sim\sum_{j=0}^\infty\frac{\delta_{2j}}{\nu^{2j}},\quad \nu \to\infty,
\end{equation}
where $\delta_0=\alpha_0=1$ and
\begin{equation}\label{eq:lagder13}
\begin{array}{ll}
\dsp{\gamma_{2j}=\phi(\zeta)\alpha_{2j}+\alpha^\prime_{2j}+\zeta \beta_{2j+1},\quad \  j=0,1,2,,\ldots,}\\[8pt]
\dsp{\delta_{2j}=\alpha_{2j}+\beta^\prime_{2j-1}+\phi(\zeta)\beta_{2j-1},\quad j=1,2,3,\ldots\,.}\\[8pt]
\end{array}
\end{equation}
Observe that $\phi(\zeta)$ depends on $\nu$, and we can rearrange the expansions and coefficients by using the relation of $\phi(\zeta)$ given in \eqref{eq:lagder11}.

In the algorithm the expression (\ref{eq:lagder10}) is used for calculating the scaled weights $\omega_i$ and from those the unscaled weights. For the 
scaled weights (Eq. (\ref{wasym2})), we need to compute $\dot{y}(z_i)$, with $x_i=z_i^2$ the zeros of $L_n^{(\alpha)}(x)$ and $y(z)$ as in  (\ref{yzl}). Considering (\ref{eq:lagder10}) this can be written
$$
\dot{y}(z_i)=(-1)^n\Frac{\nu^{\frac{\alpha}{2}+\frac{1}{12}}}{2^{\alpha-\frac12}}\left(\Frac{\bar{x}_i-1}{\zeta_i}\right)^{1/4}
\left(\nu^{-\frac23}\Ai\left(\nu^{2/3} \zeta_i\right)C(\zeta_i)
+\Ai^{\prime}\left(\nu^{2/3}\zeta_i\right)D(\zeta_i)\right),
$$
with $\tilde{x}_i=x_i /\nu$ and $\zeta_i$ the value of $\zeta$ corresponding to the node $x_i$. 
In order to calculate  $\Ai\left(\nu^{2/3} \zeta_i\right)$
and $\Ai^{\prime}\left(\nu^{2/3}\zeta_i\right)$, we take into account that 
$\nu^{2/3}\zeta_{n-i+1}=a_i+\delta_i$, where $a_i$ is 
the $i$-th negative zero of the Airy function ${\rm Ai}(x)$ and $\delta_i$ is a small correction that is 
accurately computed asymptotically with (\ref{asyze}).
We can compute the Airy function using few terms in the Taylor series
which, using ${\rm Ai}'(x)=x{\rm Ai}(x)$ and because ${\rm Ai}(ai)=0$, can be written as
\begin{equation}\label{eq:HermAiryT1}
{\rm Ai}(\nu^{2/3}\zeta_i)=\delta_i{\rm Ai}'(a_i)\left(1+\Frac{a_i}{6}\delta_i^2
+\Frac{1}{12}\delta_i^3+\Frac{a_i^2}{120}\delta_i^4+
\Frac{a_i}{120}\delta_i^5+\ldots\right).
\end{equation}
Similarly, for the derivative
\begin{equation}\label{eq:HermAiryT2}
{\rm Ai}'(\nu^{2/3}\zeta_i)={\rm Ai}'(a_i)\left(1+\Frac{a_i}{2}\delta_i^2+\Frac{1}{3}\delta_i^3
+\Frac{a_i^2}{24}\delta_i^4+
\Frac{a_i}{20}\delta_i^5+
\ldots\right).
\end{equation}

With this we only need to compute the derivative of the Airy function at the zeros of the Airy function. 
In our algorithms we use a small number of precomputed values of $\Ai'(a_i)$
(for $i=1,...,20$), and if more values are needed we compute them with the MATLAB implementation
of the Airy function with values of $a_i$ approximated by MacMahon asymptotic expansions
 for large $i$. The computation of $\Ai'(a_i)$ is well conditioned as a function of $a_i$, because $\Ai''(a_i)=0$; because of this, the
 values $\Ai'(a_i)$ can be computed with full double accuracy.
 
 \subsubsection{Bessel-type expansion}

The Bessel expansions considered in our algorithms are those described in sections 3.3 and 3.4 of \cite{Gil:2018:AAT}. The only 
aspect which was not discussed in detail in \cite{Gil:2018:AAT} is the fact that the computation of the
 Bessel functions appearing
in the expansions for the weights requires special precautions.

The evaluation of the weights involves the asymptotic computation (for large $n$) of the function $\dot{y}(z_i)$ in (\ref{wasym2}) by 
means of the expansion of \cite[Eqs. (114) and (115)]{Gil:2018:AAT}\footnote{We notice 
that we missed to mention in Eq. (115) of \cite{Gil:2018:AAT} 
that the variable $\varphi$ appearing in this equation 
is $\varphi=\frac12 (1-1/\rho)^{1/2}/\sqrt{\zeta}$.}, and this requires the computation of
the Bessel function  $J_{\alpha} \left(\nu \sqrt{\zeta}\right)$ for values of $\nu \sqrt{\zeta}$ 
near zeros of  $J_{\alpha} \left(\nu \sqrt{\zeta}\right)$\footnote{Observe that for the Bessel-type expansions in 
\cite[sections 3.3 and 3.4]{Gil:2018:AAT} we denoted $\nu=2n+\alpha+1$, while here we prefer 
to denote $\nu=4n+2\alpha+2$, as in 
Eq. (\ref{normL}), in order to 
be consistent with the notation of the Airy expansion we described earlier.}. This computation can be unstable if done directly using standard numerical implementations for
the Bessel function. 
The situation
is very similar to the one we described before for the Airy-type expansion.
 
More specifically, for computing the nodes using the Bessel expansion we use the values 
$$\zeta_k=\zeta_{0,k}+\epsilon_k=\Frac{j^2_k}{\nu^2} +\epsilon_k,$$ 
where $j_k$ is a
zero of the Bessel function $J_{\alpha}(z)$ and $\epsilon_k$ is a  small correction which
is computed accurately with an asymptotic series in powers of $\nu^{-2}$. 
The computation of the weights thus requires the computation
$J_{\alpha} \left(\nu \sqrt{\zeta_k}\right)$. 
We write
\begin{equation}
\nu \sqrt{\zeta_k}=j_k+\delta_k,\, \delta_k =\Frac{\epsilon_k/\zeta_{0,k}}{1+ \sqrt{1 +\epsilon_k/\zeta_{0,k}  }}j_k,
\end{equation}
and $\delta_k$ is accurately computed. Then, for computing $J_{\alpha}(j_k+\delta_k)$, where $J_{\alpha}(j_k)=0$ and 
$\delta_k/j_k<<1$ we proceed similarly as in \cite[section 3.3]{Gil:2019:NCO}, using the expansion 

\begin{equation}
J_\alpha(u+h)=\lambda^\alpha\sum_{m=0}^\infty \frac{w^m}{m!}J_{\alpha+m}(u),\quad w=-\frac{h(2u+h)}{2u}.
\end{equation}
with $u=j_k$ and $h=\delta_k$. 
 
\subsection{Combining the iterative and asymptotic methods}

For deciding the regions where each of the two methods is used in the final algorithm we follow a similar strategy as for the 
Gauss--Jacobi case. 

We compare both methods against a higher precision algorithm. Comparing
against an extended precision of the Golub-Welsch algorithm is not convenient because the exponential decay of the weights causes 
drastic loss of the accuracy for the smallest weights\footnote{For example, for
 $n=30$, only the first $20$ weights (corresponding to the smallest nodes) are computed with relative accuracy better
 than $10^{-10}$ in double precision arithmetic.}.
In addition, 
scaling of the weights is not possible in the GW method. Instead, we
 compare the asymptotic and iterative methods against a quadruple precision version of the iterative method.  In these tests we 
 consider subsampling and only those nodes $x_i$ and
  weights $w_i$ are computed such that, approximately, $w_{i}/\max_{j=1,\ldots n}\{w_j\}>10^{-300}$ 
  (more precisely, such that
  $F_i>-300\log(10)$, see (\ref{efe})).

\begin{figure}[htbp]
\begin{center}
\includegraphics[width=0.49\textwidth]{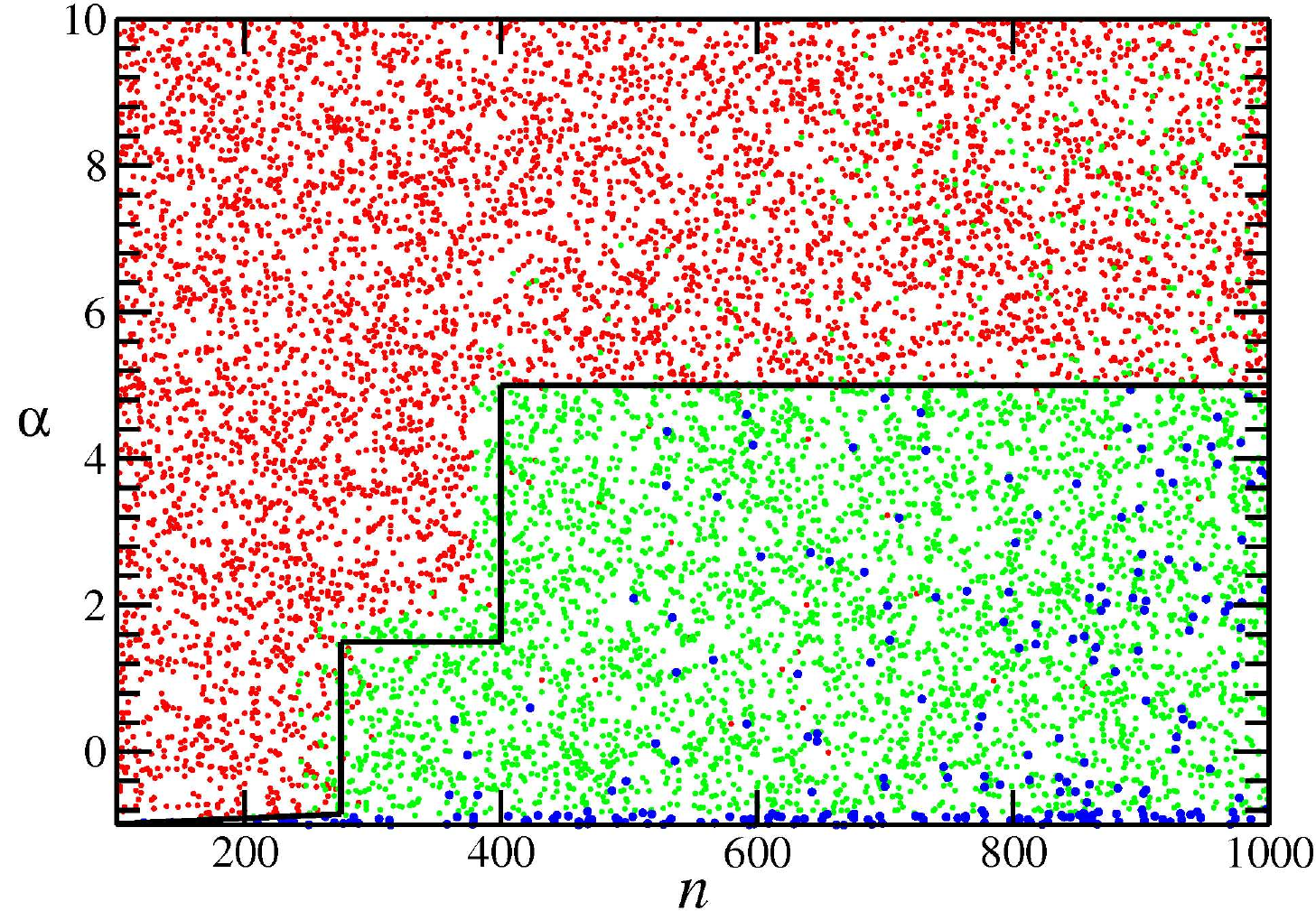}
\includegraphics[width=0.49\textwidth]{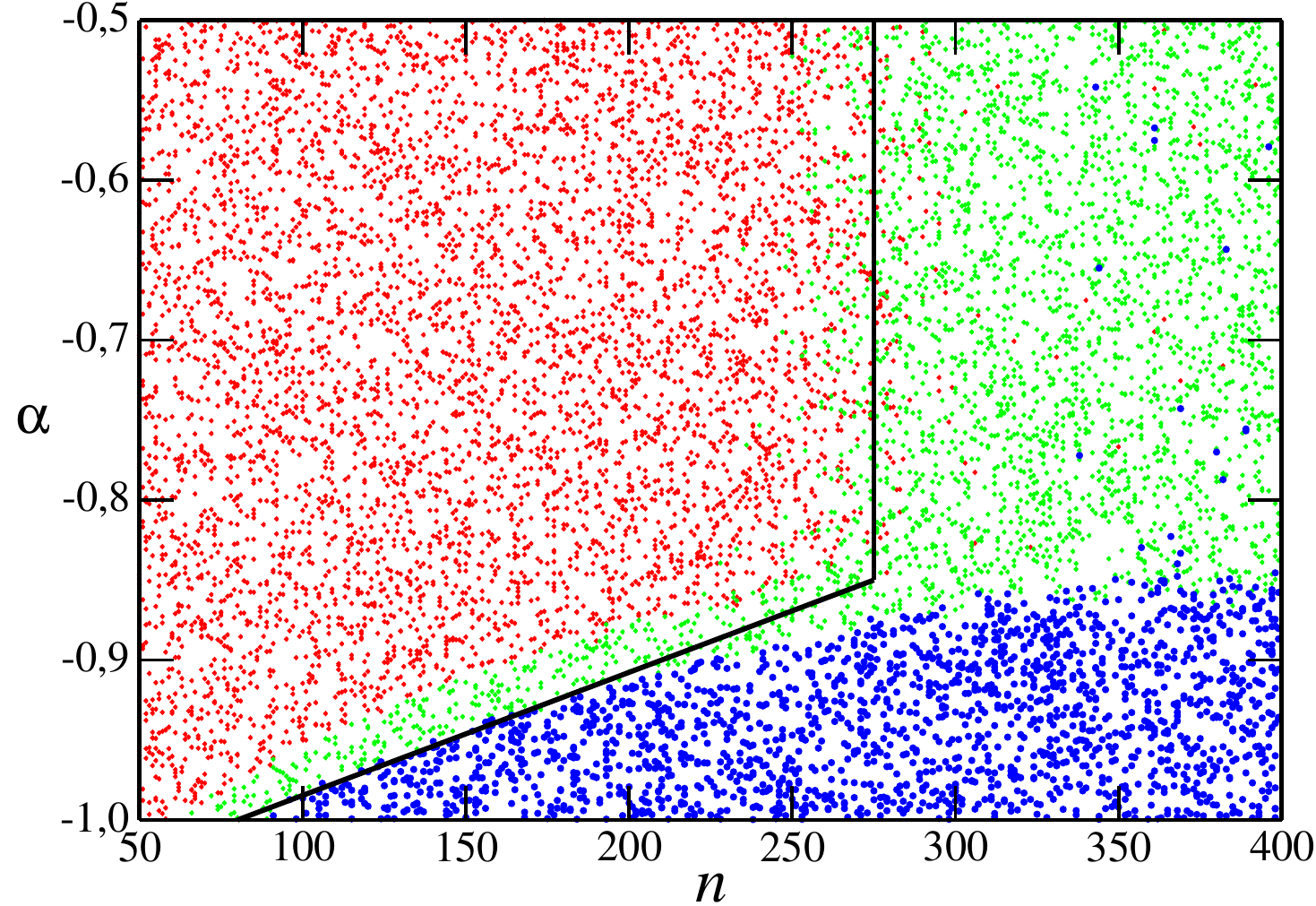}
\end{center}
\caption{
Comparison of the maximum relative error of the iterative and the asymptotic method for computing Gauss--Laguerre quadratures.
 At the red points 
 the iterative method is more
 accurate that the asymptotic method, while the opposite is true both for the green and the blue points; at the blue points 
 the accuracy of the iterative method is ten times worse than that of the asymptotic method. The separation lines between the methods
 used in the combined algorithm are also shown.
}
\label{regioL}
\end{figure}

  Fig. \ref{regioL} shows the regions where the iterative method is more accurate than the asymptotic method and vice-versa.
  As measure of the error, we take the maximum
  relative error for the scaled weights (with unitary normalization); the other two error measures considered before in the 
  Gauss--Jacobi case (average relative error and relative total error) give similar pictures but favoring slightly more the asymptotic
  approach. We give preference to the iterative method because it is a faster and so 
  we consider the error criteria favoring this method.

From the comparison of accuracies of Fig. \ref{regioL}, it appears reasonable to apply
 the asymptotic method in the following regions of parameter space:
\begin{enumerate}
\item{}$80<n\le 275$ and $-1<a\le (n-1380)/1300$.
\item{}$275<n\le 400$ and $-1<a\le 1.5$.
\item{}$400<n$ and $-1<a\le 5$.
\end{enumerate}

For the rest of parameters we use the iterative method except $n\le 5$, where we use the GW method. 
We do not consider the iterative method for $n<5$ because the method was not originally designed for such small values of $n$ (see   
\cite[section 4.2 and 4.3 ]{Gil:2019:FRA} for further details). A version for $n<5$ is possible, but GW does the job perfectly for
any $\alpha$ 
(no scaled weights are produced, though).

The asymptotic expansions are computed with a fixed number of terms. The performance may vary by adding additional terms, particularly
in the small $n$ region, at the cost of less computational efficiency. Because the iterative method is quite accurate and more 
efficient that the asymptotic method (roughly $2-3$ times faster), we
prefer to use the latter method when possible. We notice, however, that the asymptotic method is more accurate for moderate $\alpha$, 
and particularly as we approach the boundary $\alpha=-1$. This is as expected, because in the limit $\alpha\rightarrow -1^+$, the
smallest root of the Laguerre polynomial approaches zero, which is a singularity of the ODE; therefore, the Taylor method used for
computing the solution of the ODE (\ref{yzl}) in the iterative method tends to fail. This is clearly shown in Fig. \ref{regioL} (right)
where it is observed that the asymptotic method is better for small $\alpha$ even for values close to $n=80$.

As commented in the Gauss--Jacobi case, 
the accuracy of the Bessel-type asymptotic expansions is limited by the accuracy
in the computation of Bessel functions with the intrinsic MATLAB procedure. 
For the Gauss--Laguerre case, this limitation is of higher importance because it affects
the smallest positive nodes, and therefore the most significant weights for small $\alpha$.
In order to improve the accuracy
of the algorithm we use our own algorithms (included in the package) for the computation of Bessel functions, which prove to be
more accurate that the MATLAB intrinsic implementation, at least for small orders; we use our expansions for $-1<\alpha \le 5$ and use
the intrinsic MATLAB codes for $\alpha>5$.

\subsection{Testing the combined algorithm}

\begin{figure}[htbp]
\begin{center}
\includegraphics[width=0.49\textwidth]{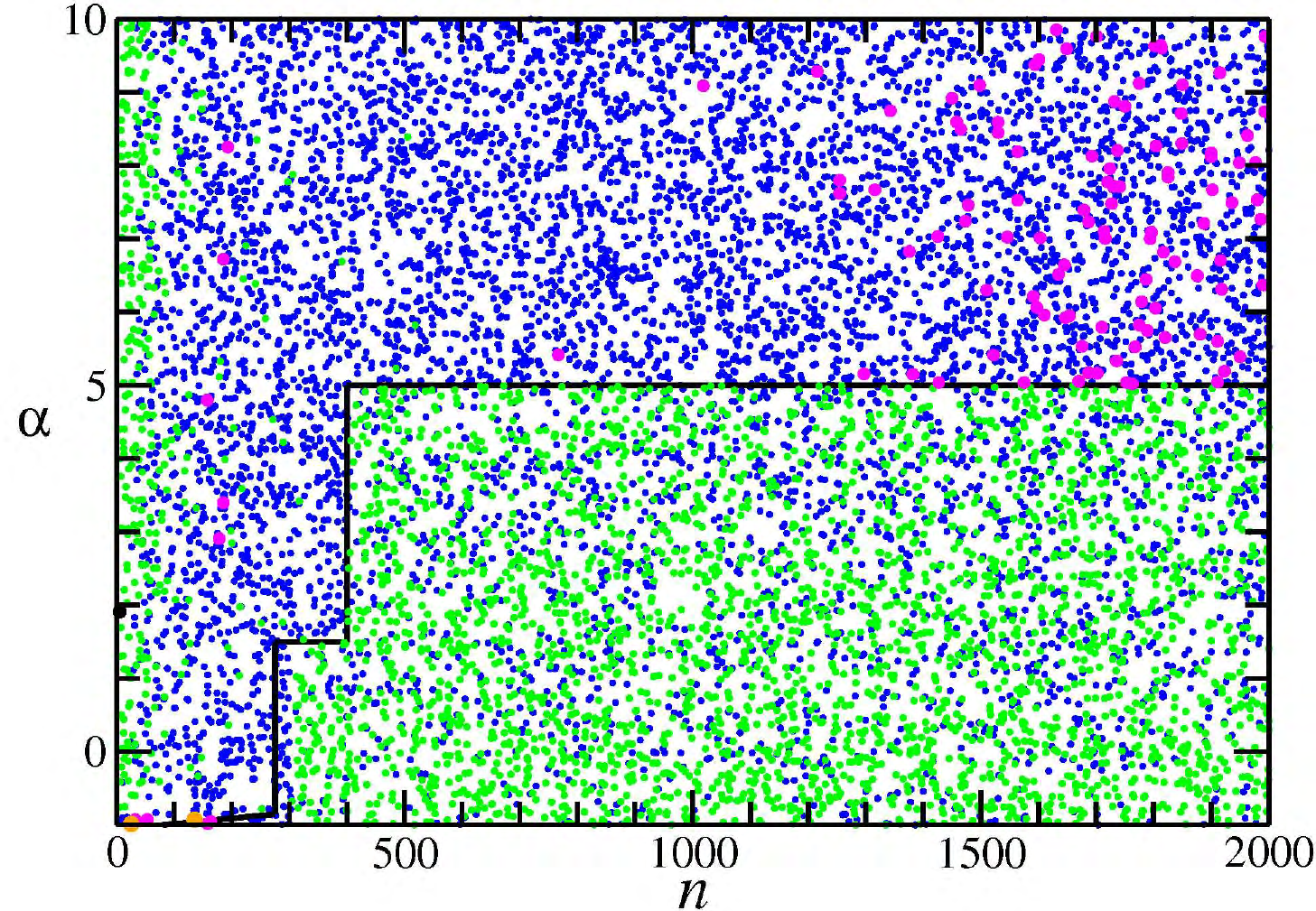}
\includegraphics[width=0.49\textwidth]{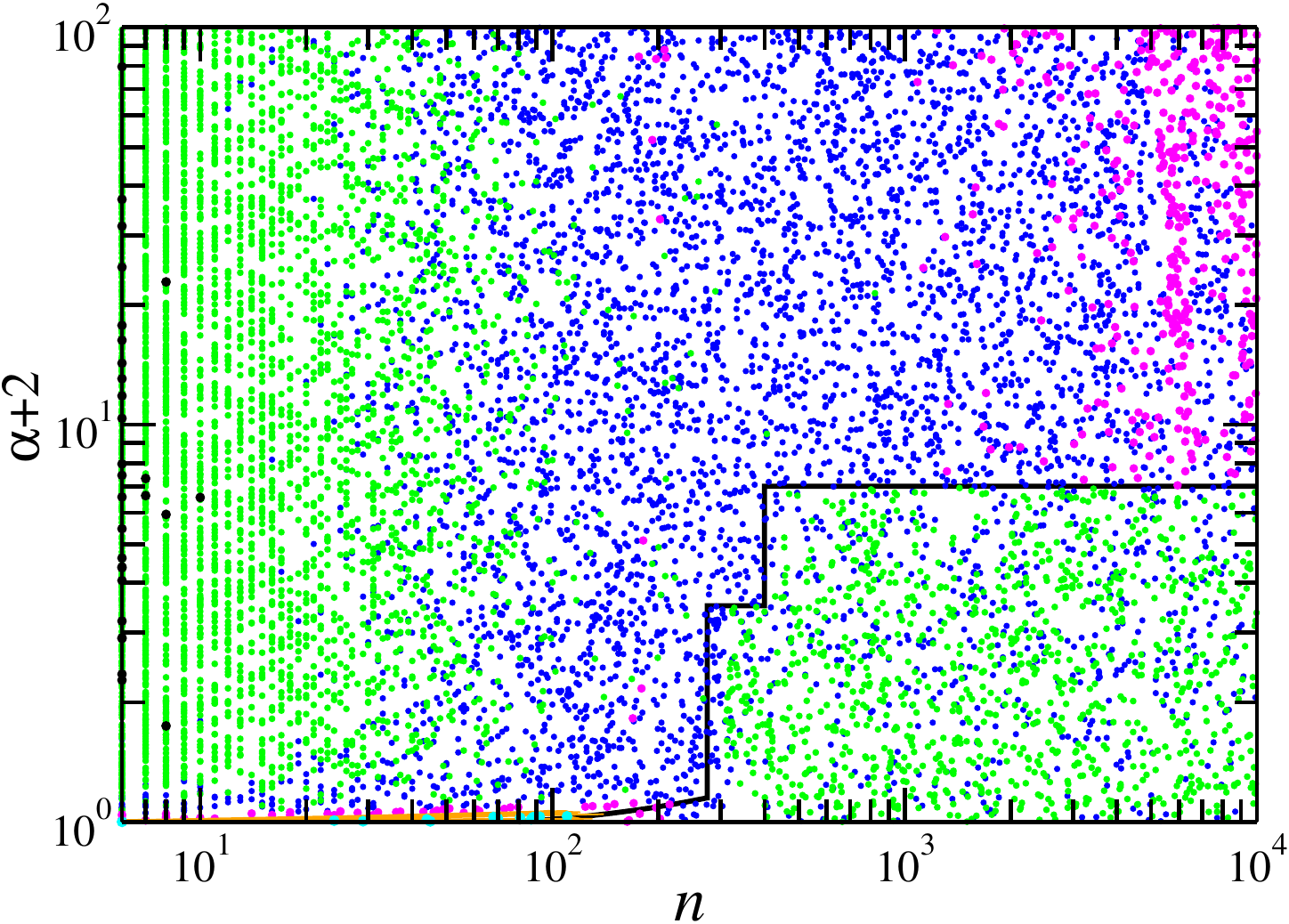}
\end{center}
\caption{
Accuracy of the scaled weights obtained with the combined algorithm for two ranges of the parameters. The different colors 
represent maximum relative accuracies in the intervals $[10^{-16+k},10^{-16+k+1})$, with $k=0$ for back, $k=1$ for green,
$k=2$ for blue, $k=3$ for magenta, $k=4$ for orange and $k=5$ for cyan. The figures also show the separation line between the
asymptotic and iterative methods.
}
\label{accL}
\end{figure}

Fig. \ref{accL} shows the accuracy of the combined algorithm for computing scaled unitary weights when subsampling is considered
(with the same level of subsampling as in the previous section). The figure on the left is generated with random points uniformly distributed, while the figure on the right shows a larger range with points exponentially generated. Both figures show that the asymptotic
method produces accuracies for the scaled weights 
always better than $10^{-13}$ and most of the times better than $10^{-14}$. The worst case scenario for the algorithm is when 
the iterative method is used for values very close to $\alpha=-1$, for the reasons explained before.

For brevity, we do not show a detailed comparison with other methods, and in particular with the algorithm {\bf lagpts} included in the Chebfun 
package \cite{Driscoll:2013:chebfun}. Our algorithm is generally faster and also more accurate. 
More specifically, for $n<3000$, where {\bf lagpts} used a combination of non-asymptotic methods, 
our codes are typically $5-10$ times faster and
compute the nodes and weights accurately in a much larger range of $\alpha$ (practically unrestricted). 
For $n>3000$ {\bf lagpts} uses the asymptotics approximations from \cite{Opsomer:2023:HOA}
and it might be faster (by less than a factor $2$) than our algorithm when asymptotics is used\footnote{For such large values of $n$ our algorithm may be accelerated by considering
less terms in the asymptotic expansions; we will consider this possibility in future versions of the software}, 
but we notice a significant drop in the accuracy of the weights given by {\bf lagpts}. 
For instance, for $n=3002$, $\alpha=0$, the worst relative error for the first $100$ weights
(the most significant ones) is $3\,10^{-15}$ for our algorithm and $9\,10^{-10}$ for {\bf lagpts}; the situation gets worse as 
the degree becomes larger.
On the other hand, our algorithm presents novel features not shared by previous algorithms as there are no practical restrictions on the
parameters and subsampling is available. 
 
\section{Gauss--Hermite quadrature}

The algorithms in this section are based on the asymptotic methods of 
\cite{Gil:2018:AAT} and the iterative methods of \cite{Gil:2019:FRA}. 
We outline the main ingredients of the iterative method needed and
 we describe further details of the asymptotic method that were not fully developed in \cite{Gil:2018:AAT}.
 
 Hermite polynomials $H_n (x)$ are solutions of the differential equation $w''(x)-2xw'(x)+2nw(x)=0$. With the 
 change of function $y(x)=e^{-x^2 /2}H_n(x)$, the transformed function satisfies
 \begin{equation}
 \label{ODEH}
 y''(x)+A(x)y(x)=0,\,A(x)=2n+1-x^2.
 \end{equation}
The nodes $x_i$ of $n$-point Gauss--Hermite quadrature are the roots of $H_n(x)$, which are symmetric around the origin, and the weights $w_i$
can be written as (see Table \ref{table1}):
$$
w_i =e^{-x_i^2}\omega_i,\,\omega_i=\sqrt{\pi}2^{n+1}n!y'(x_i)^{-2},
$$
where the quantities $\omega_i$ are the scaled weights. The scaled weights have slow variation as a function of $i$ 
compared with the unscaled weights $w_i$, in accordance with the asymptotic estimate for the weights
$$
w_i\sim \Frac{\pi}{\sqrt{2n}}e^{-x_i^2},\,n\rightarrow \infty,
$$
and as could be expected from the behavior of the ODE (\ref{ODEH}).

\subsection{Iterative method}

The iterative method for Gauss--Hermite is a simple, accurate and efficient method of computation. 
The nodes are computed with full double accuracy, while
some mild degradation is observed in the computation of the weights as the 
degree increases. For computing scaled weights for large $n$, the asymptotic methods 
produce higher accuracy. However, as we discuss later, 
the iterative method for Gauss--Hermite also has a number of advantages with respect to asymptotics.

The iterative method is based on the use of fixed point method stemming from the ODE (\ref{ODEH}), and the positive nodes
 are generated in increasing order together with their corresponding weights. The algorithm is
 remarkably simple and accurate and can be prepared to work with arbitrary accuracy. We provide both MATLAB (double precision) 
 and Maple (arbitrary precision) implementations. For a detailed description of the algorithm we refer to \cite{Gil:2019:FRA}.   

\subsection{Asymptotic methods}

We describe here some numerical aspects of the asymptotic computation which
were not described in \cite{Gil:2018:AAT} and which are important for an
accurate computation, particularly of the weights.

We consider two types of asymptotic methods: asymptotics in terms of elementary 
functions and Airy-type asymptotics. The former is accurate for most of the nodes
while the accuracy degrades for the largest nodes, for which the Airy-type 
expansions are preferable. We have checked numerically that only the last $20$ 
nodes (at most) should be computed with the Airy-type expansion.

The use of the elementary expansion is described in detail in \cite{Gil:2019:FRA}; 
however, for the Airy-type expansion some improvements are considered in the algorithm
which we describe next.

\subsubsection{Computation of the Airy-type expansion}

As in \cite{Gil:2018:AAT}, we denote $y(x)=e^{-x^2/2}H_n (x)$ and define the scaled weights as
$$
\tilde{w}_i=\Frac{\sqrt{\pi}2^{n+1}n!}{y'(x_i)^2}=2\sqrt{\pi}n! \left[U'\left(-n-\frac12,\sqrt{2}x_i\right)\right]^{-2}
$$
where in the last equation we have used the relation between Hermite polynomials and 
parabolic cylinder functions
\cite[12.7.2]{Temme:2010:PCF} $y(x)=2^{n/2}U(-n-\frac12,\sqrt{2}x)$. As we did earlier in the Laguerre case, we need to give
a brief description regarding the asymptotic computation of the derivative $y'(x_i)$. For this purpose, we use 
the Airy-type expansion for Parabolic Cylinder functions (see \cite{Temme:2010:PCF}, section 12.10(vii)). 

Denoting $\mu=\sqrt{2n+1}$ and $t=x/\mu$, we have 
\begin{equation}\label{eq:Hp05}
U'\left(-n-\frac12,\sqrt{2}x\right)\sim\frac{(2\pi)^{\frac{1}{2}}%
\mu^{\frac{2}{3}}g(\mu)}{\phi(\zeta)}\*\left(\frac{\operatorname{Ai}\left(\mu^%
{\frac{4}{3}}\zeta\right)}{\mu^{\frac{4}{3}}}C(\zeta)+\operatorname{Ai}'\left(\mu^{\frac{4}{3}}\zeta\right)D(\zeta)\right),
\end{equation}
where
$
\frac13 (-\zeta)^{3/2}=\Frac{1}{2}\left(\arccos t -t\sqrt{1-t^2}\right)
$
in the region where the nodes lie ($|t|<1$) and $g(\mu)$, $\phi(\zeta)$ are as defined in \cite[section 12.20]{Temme:2010:PCF}. The
coefficients admit the expansions
\begin{equation}\label{eq:Hp06}
\begin{array}{ll}
&C(\zeta)=\dsp{\sum_{s=0}^{\infty}\frac{C_{s}(\zeta)}{\mu^{4s}},\quad D(\zeta)=
\sum_{s=0}^{\infty}\frac{D_{s}(\zeta)}{\mu^{4s}}.}\\[8pt]
\end{array}
\end{equation}

In order to avoid numerical cancellations when computing the coefficients when $t=x/\mu$ 
is close to $1$ ($\zeta$ small) in the above representations, which is the case for the largest zeros,
we can first expand $\zeta$ in powers of  $\vert t-1\vert$ and the
coefficients  
$C_{k}(\zeta),\,D_{k}(\zeta)$, which are analytic at $\zeta=0$, in powers of $\wt\zeta=2^{-\frac13}\zeta$.

For small values of $\vert t-1\vert$ and $\vert\zeta\vert$, we use the expansions
\begin{equation}\label{eq:HermAiry06}
\begin{array}{l}
2^{-\frac{1}{3}}\zeta=\dsp{(t-1)+\frac{1}{10}(t-1)^2-\frac{2}{175}(t-1)^3+{O}\left((t-1)^4\right)}\\
\\
C_k(\zeta)=2^{-\frac13}\dsp\sum_{j=0}c^k_j{\wt\zeta}^j,\,\, D_k(\zeta)=\dsp\sum_{j=0}d^k_j{\wt\zeta}^j,  
\end{array}
\end{equation}
where $\wt\zeta=2^{-\frac{1}{3}}\zeta$ and the coefficients $c_j^k$, $d_j^k$ are rational numbers.

With respect to the computation of the Airy functions, a similar procedure as described for the Gauss--Laguerre case 
can be used here. Indeed, the dominant term in the asymptotic approximation for the 
nodes is determined by setting to zero the dominant contribution in the expansion for the Hermite polynomial 
\cite[12.10.35]{Temme:2010:PCF}.
If $\zeta_i$ is the value of $\zeta$ corresponding to the node $x_{i}$, we have
$
\mu^{4/3}\zeta_{n-i+1}=a_i +\delta_i
$,
where $\delta_i$ can be computed accurately and $\delta_i/a_i<<1$. We can compute the Airy functions in (\ref{eq:Hp05}) using Taylor series. For Gauss--Hermite quadrature the Airy-type expansion is
only used for the last $20$ weights and we store the precomputed values of $a_i$ and ${\rm Ai}'(ai)$, $i=1,\ldots 20$.

\subsection{Comparing the iterative and asymptotic methods}

The accuracy of the iterative and asymptotic methods is compared in Table \ref{tablerr}. 
The nodes are computed with slightly more accuracy for the iterative method, which produces full double precision
accuracy for any degree. Contrarily, the scaled weights
are more accurately computed with the asymptotic method when $n\ge 150$ than with the iterative method. 
However, the maximum relative error for the unscaled weights which are
computable in double precision (limited by underflow) is very similar for the iterative and asymptotic method. 
The main source of error for the unscaled weights is the exponential factor $e^{-x_i^2}$, which is slightly more accurately 
computed with the iterative
method. 

Therefore, each
method has its advantages in terms of accuracy. But because the iterative method is faster by a factor $10$ approximately, we consider
the iterative method the main method of computation. Differently to the Jacobi and Laguerre cases, we do not combined the two approaches in 
a single algorithm because none of them is superior in terms of accuracy to the other one, but the iterative method is faster.

\begin{table}
$$
\begin{array}{|c|c||c|c||c|c|c|}
\hline
n    & \max\epsilon (x_i)  & \max\epsilon (\omega_i) & \epsilon (w_1) & \% & \max\epsilon (\omega_i) & \max \epsilon (w_i)   \\ 
\hline  
150  &  1.6\,10^{-16}  &  1.5\,10^{-14} & 3.2\,10^{-15}  & 100   & 1.5\,10^{-14} & 1.1\,10^{-13}  \\ 
     &  5.8\,10^{-16}  &  2.2\,10^{-15} & 5.5\,10^{-16}  &       & 2.2\,10^{-15} & 9.5\,10^{-14}   \\
\hline     
500  &  1.7\,10^{-16}  &  1.5\,10^{-14} & 2.1\,10^{-15}  & 92    & 1.5\,10^{-14} & 1.2\,10^{-13}   \\  
     &  5.8\,10^{-16}  &  1.3\,10^{-15}  & 3.7\,10^{-16} &       & 1.3\,10^{-15} & 3.2\,10^{-13}  \\
\hline     
10^3 &  1.7\,10^{-16}  &  3.8\,10^{-14}  & 3.1\,10^{-15} & 70.2  & 2.7\,10^{-14} & 1.3\,10^{-13}  \\
     &  7.7\,10^{-16}  &  2.9\,10^{-15} & 7.5\,10^{-17}  &       & 9.4\,10^{-16} & 6.8\,10^{-13}  \\
\hline     
10^4 &  2.0\,10^{-16}  &  7.7\,10^{-13} & 3.0\,10^{-15}  & 23.52 & 1.3\,10^{-13} & 2.3\,10^{-13}  \\
     &  1.2\,10^{-15}  &  3.8\,10^{-15}  & 1.1\,10^{-15} &       & 1.8\,10^{-15} & 8.8\,10^{-13}  \\
\hline     
10^5 &  1.6\,10^{-16}  &  8.4\,10^{-13}  & 2.0\,10^{-14} & 7.47  & 2.6\,10^{-13} & 3.9\,10^{-13}  \\
     &  1.6\,10^{-15}  &  4.7\,10^{-15}  &  7.1\,10^{-16} &       & 1.3\,10^{-15} & 4.9\,10^{-13} \\
\hline     
10^6 &  1.6\,10^{-16}  &  5.1\,10^{-11}  & 6.4 \,10^{-14} & 2.36  & 7.1\,10^{-13} & 8.3\,10^{-13}  \\
     &  1.5\,10^{-15}  &  4.4\,10^{-15}  &  7.5\,10^{-16} &       & 1.4\,10^{-15} & 6.8\,10^{-13} \\     
\hline
\end{array}
$$
\caption{Relative errors for the iterative and asymptotic methods for different values of the degree $n$. For each $n$, the
first line corresponds to the iterative method and the second to asymptotic method. $\max\epsilon (x_i)$,  
$\max\epsilon (\omega_i)$ and $\max\epsilon (w_i)$ denote the maximum relative error 
for the nodes, scaled weights and unscaled weights respectively, while 
$\epsilon(w_1)$ is the relative error for the 
most significant weight. The last three columns correspond to the computation with subsampling, in which only the nodes and weights 
such that 
$w_i/w_1>10^{-300}$ are computed; the column labeled as $\%$ gives the percentage of nodes that are computed.}
\label{tablerr}
\end{table}

\section{Further applications}

In this section, we discuss how to compute the Gauss--Lobatto and Gauss--Radau quadratures, as well as the
barycentric weights associated to the different sets of classical nodes (including Gauss--Lobatto and Gauss--Radau variants).

\subsection{Gauss--Radau and Gauss--Lobatto quadratures}

In some applications, for instance for solving boundary problems, it can be interesting to consider rules with the 
maximal degree of polynomial exactness which include one or two of the endpoints
of the interval of integration as nodes of integration. When only one of the endpoints is considered the rule is called of
Gauss--Radau type, and when the two endpoints are used of Gauss--Lobatto type. 

Let us first consider a generic Gauss--Radau case for a weight function $w(x)$ in an interval $[a,b]$ with $x=a$ as
a node:
\begin{equation}
\label{GR}
I(f)=\displaystyle\int_a^b f(x) w(x) dx\approx Q_n (f)=w_0 f(a) + \displaystyle\sum_{i=1}^n w_i f(x_i).
\end{equation}
The maximal possible degree of accuracy is $2n$. We describe briefly how to obtain the nodes and weights in this
generic case and we later summarize the results for the classical cases.

According to a well-know theorem (see the first theorem in 
section 2.7.1 of \cite{Davis:1984:MON})
 the maximal degree in rule (\ref{GR}) is achieved with an interpolatory formula which has
as internal nodes those of the (for instance monic) polynomial $p_n^R(x)$ of degree $n$ which is orthogonal to all the polynomials 
of degree less than $n$ with respect to the weight function $w^{R} (x)=(x-a)w(x)$.

Notice that, denoting by $w_i^{R}$ the weights of the $n$-point Gaussian quadrature for the weight function $w^R (x)$ in 
the interval $[a,b]$, we have
$$
I_k=\displaystyle\int_a^b x^k w^R (x)dx=\displaystyle\sum_{i=1}^n w_i^{R} x_i^k,\,k=0,1,\ldots 2n-1.
$$
On the other hand, because $w^{R} (x)=(x-a)w(x)$ and the rule $Q_n (f)$ has degree of exactness $2n$
$$
I_k=I((x-a)x^k)=\displaystyle\sum_{i=1}^n w_i (x_i-a) x_i^k,\,k=0,1,\ldots 2n-1.
$$
Therefore, considering the last two equations one concludes that 

\begin{equation}
\label{nodesR}
w_i=w_i^{R}/(x_i-a).
\end{equation}

For the Gauss--Radau case with $x=b$ as a node we proceed similarly. For the Gauss--Lobatto case, with both $x=a$ and $x_b$
as nodes, the internal nodes are the zeros of the polynomial $p_n^{L}(x)$ 
orthogonal with respect $w^{L}(x)=(x-a)(b-x)w(x)$, and the relation between
weights follow similarly: $w_i=w_i^{L}/[(x_i-a)(b-x_i)]$.

As for the boundary node(s), and considering again the Gauss--Radau case with $x=a$ as a node, because $p^R_n(x_i)=0,\,i=1..n$ and
$Q_n$ has degree of exactness $2n$
\begin{equation}
\label{bou}
I(p_n^R (x))=\displaystyle\int_a^b p^R_n(x) w(x) dx = Q_n (p_n^R (x))=w_0 p_n^R(a).
\end{equation}

For the Gauss--Lobatto case we would consider the integration of $(x-a)p_n^L (x)$ and  $(b-x)p_n^L (x)$ 
for determining the boundary
nodes.

For the cases of classical Gaussian rules the polynomials $p_n^R (x)$ are classical orthogonal
 polynomials themselves, and the boundary nodes can be computed explicitly. 
  Take for instance the case of Gauss--Radau--Laguerre quadrature, 
 which is the Radau case considered before with $a=0$, $b=+\infty$, $w(x)=x^{-\alpha} e^{-x}$. 
 The Laguerre weight is $w(x)=x^\alpha e^{-x}$, and 
 in this case $p_n^R(x)=L_{n}^{(\alpha+1)}(x)$ is the orthogonal polynomial associated to the weight $w^R (x)=x w(x)$. To compute
 the Gauss--Radau--Laguerre quadrature, we first compute the nodes $x_i^R$ and weights 
 $w_i^R$ of Gauss--Laguerre quadrature with 
 parameter $\alpha+1$; then the internal nodes of the Gauss--Radau--Laguerre quadrature with parameter $\alpha$ are
 $x_i=x_i^R$, $i=1,\ldots n$, and the internal weights are $w_i=w_i^{R}/x_i$, $i=1,\ldots n$. The boundary weight, considering 
 (\ref{bou}) is
 \begin{equation}
 \label{garal}
 w_0=\Frac{1}{L_n^{(\alpha+1)}(0)} \displaystyle\int_0^{+\infty}L_n^{(\alpha+1)}(x)x^\alpha e^{-x}dx,
 \end{equation}
 and because $L_n^{(\alpha+1)}(x)=\sum_{j=0}^{n}L_{j}^{(\alpha)}(x)$ \cite[18.18.37]{Koornwinder:2010:OP}, using
 the orthogonality of Laguerre polynomials and the fact that $L_n^{(\alpha)}(0)=(\alpha+1)_n/n!$ we have
 \begin{equation}
 \label{w0l}
 w_0=\Frac{n!}{(\alpha+2)_n}\displaystyle\int_0^{+\infty}x^\alpha e^{-x}dx=\Frac{\Gamma(\alpha+1)
 \Gamma(\alpha+2)n!}{\Gamma(n+\alpha+2)}.
 \end{equation}
 
For the rest of Gauss--Radau and Gauss--Lobatto quadratures the computation of the boundary weights follow similar ideas, and the
main task is to compute, as in (\ref{garal}), the integral of an orthogonal polynomial of the same family associated to the
weight function, but with one displaced parameter. For details, we mention 
\cite{Gautschi:2000:GRF} for Gauss--Radau formulae and \cite{Gautschi:2000:HOG} for Gauss--Lobatto, and we only summarize the results:

1. Gauss--Laguerre--Radau.
Let $x_i^R$ and $w_i^R$ be the nodes and weights for the $n$-point Gauss--Laguerre formula with parameter $\alpha+1$, then
the quadrature rule
\begin{equation}
\label{GLR}
\displaystyle\int_0^{+\infty} f(x)x^{\alpha}e^{-x}dx\approx Q_{n}(f)=w_0 f(0) +\displaystyle\sum_{i=1}^n w_i f(x_i)
\end{equation}
has maximal degree of exactness $2n$ taking $x_i=x_i^R$, $w_i=w_i^R/x_i^R$, $i=1,\ldots n$, and with 
$w_0$ given by (\ref{w0l}).

2. Gauss--Jacobi--Radau with $x=-1$ as a node.
Let $x_i^R$ and $w_i^R$ be the nodes and weights for the $n$-point Gauss--Jacobi formula with parameters $\alpha$ and 
$\beta+1$, then
the quadrature rule
\begin{equation}
\label{GJRL}
\displaystyle\int_{-1}^{1} f(x)(1-x)^\alpha (1+x)^\beta dx\approx Q_{n}(f)=w_0 f(-1) +\displaystyle\sum_{i=1}^n w_i f(x_i)
\end{equation}
has maximal degree of exactness $2n$ taking $x_i=x_i^R$, $w_i=w_i^R/(1+x_i^R)$, $i=1,\ldots n$, and with 
$$
w_0=2^{\alpha+\beta+1}\Frac{\Gamma (\beta+1)\Gamma (\beta+2)\Gamma(n+1)\Gamma (n+\alpha+1)}
{\Gamma (n+\beta+2)\Gamma(n+\alpha+\beta+2)}.
$$

3. Gauss--Jacobi--Radau with $x=+1$ as a node.
Let $x_i^R$ and $w_i^R$ be the nodes and weights for the $n$-point Gauss--Jacobi formula with parameters $\alpha +1$ and 
$\beta$, then
the quadrature rule
\begin{equation}
\label{GJRR}
\displaystyle\int_{-1}^{1} f(x)(1-x)^\alpha (1+x)^\beta dx\approx Q_{n}(f)=\displaystyle\sum_{i=1}^n w_i f(x_i)\,+w_{n+1}f(1)
\end{equation}
has maximal degree of exactness $2n$ taking $x_i=x_i^R$, $w_i=w_i^R/(1-x_i^R)$, $i=1,\ldots n$, and with 
$$
w_{n+1}=2^{\alpha+\beta+1}\Frac{\Gamma (\alpha+1)\Gamma (\alpha+2)\Gamma(n+1)\Gamma (n+\beta+1)}
{\Gamma (n+\alpha+2)\Gamma(n+\alpha+\beta+2)}
$$ 

4. Gauss--Jacobi--Lobatto.
Let $x_i^R$ and $w_i^R$ be the nodes and weights for the $n$-point Gauss--Jacobi formula with parameters $\alpha +1$ and 
$\beta +1$, then
the quadrature rule
\begin{equation}
\label{GJL}
\displaystyle\int_{-1}^{1} f(x)(1-x)^\alpha (1+x)^\beta dx\approx Q_{n}(f)=w_0 f(-1)+
\displaystyle\sum_{i=1}^n w_i f(x_i)\,+w_{n+1}f(1)
\end{equation}
has maximal degree of exactness $2n+1$ taking $x_i=x_i^R$, $w_i=w_i^R/(1-(x_i^R)^2)$, $i=1,\ldots n$, and with 
$$
w_{0}=2^{\alpha+\beta+1}\Frac{\Gamma (\beta+1)\Gamma (\beta+2)\Gamma(n+1)\Gamma (n+\alpha+2)}
{\Gamma (n+\beta+2)\Gamma(n+\alpha+\beta+3)},
$$ 
$$
w_{n+1}=2^{\alpha+\beta+1}\Frac{\Gamma (\alpha+1)\Gamma (\alpha+2)\Gamma(n+1)\Gamma (n+\beta+2)}
{\Gamma (n+\alpha+2)\Gamma(n+\alpha+\beta+3)}.
$$

\subsection{Barycentric interpolation}

Let $S=\left\{(x_i,f_i),\,i=1,\dots n\right\}$ be a set of points of ${\mathbb R}^2$ with $x_i\neq x_j$ if $i\neq j$, then 
it is well known that the polynomial $P(x)$ of the smallest degree such that $P(x_i)=f_i$ can be written in the barycentric form
as
\begin{equation}
\label{bari}
P(x)=\left.\displaystyle\sum_{i=1}^n \frac{u_i f_i}{x-x_i}\right/\displaystyle\sum_{i=1}^n\frac{u_i}{x-x_i},
\, u_i=k\left(\prod_{j=1,j\neq i}^n
(x_i-x_j)\right)^{-1},
\end{equation}
and the values $u_i$ are called barycentric weights, with $k$ an arbitrary constant. The barycentric weights can also be written as 
\begin{equation}
\label{bari2}
u_i=k \Phi' (x_i)^{-1},\,\Phi (x)=\displaystyle\prod_{j=1}^n(x-x_j).
\end{equation}
Therefore, for interpolation at the $n$ simple zeros of a polynomial of degree $n$, 
the barycentric weights can be computed as the inverse of the derivative of the polynomial 
at the zeros.

In our algorithms, we compute barycentric nodes and weights for the classical Gaussian quadratures, including also the 
Gauss--Radau and Gauss--Lobatto variants. 
Just by comparing the expression of the Gaussian weights with (\ref{bari2}) we can relate the barycentric weights with the Gaussian
weights. $\Phi'(x_i)$ in (\ref{bari2}) is the derivative of the monic orthogonal polynomial at the node $x_i$, 
and this same quantity appears in the expression for the Gauss weights (squared). Neglecting constant multiplicative factors 
(not depending on the nodes $x_i$), 
the relation is thus very simple.   For the Hermite case we can write 
$u_i=(-1)^i \sqrt{w_i}$ where the $(-1)^i$ is due to the fact that, being all the zeros of $\Phi(x)$ simple, the sign of the derivative
is opposite in consecutive zeros; similarly, for the Laguerre case $u_i=(-1)^i\sqrt{x_iw_i}$ and for the Jacobi case 
$u_i=(-1)^i\sqrt{(1-x_i^2)w_i}$.

It is also possible to relate the Gauss and barycentric weights for Gauss--Radau and Gauss--Laguerre quadratures 
(see \cite{Wang:2014:EBW}). We summarize these results in Table 3.

\begin{table}
\begin{center}
\begin{tabular}{|c|c|c|}
\hline
 Hermite & Laguerre &  Jacobi\\
 \hline
  $u_i=(-1)^i \sqrt{w_i}$ & $u_i=(-1)^i \sqrt{x_i w_i}$ & $u_i=(-1)^i\sqrt{(1-x_i^2)w_i}$\\ 
 \hline
  & $u_0=\sqrt{(\alpha+1)w_0}$ & $u_0=\sqrt{2(\beta+1)w_0}$ \\
 & $u_i=(-1)^i \sqrt{w_i}$  & $u_i=(-1)^i \sqrt{(1-x_i)w_i}$ \\
\hline
&  & $u_i=(-1)^i \sqrt{(1+x_i)w_i}$\\
 & &  $u_{n+1}=(-1)^{n+1}\sqrt{2(\alpha+1)w_{n+1}}$ \\
\hline
   &  &  $u_0=\sqrt{(\beta+1)w_0}$ \\
 & &  $u_i=(-1)^i \sqrt{w_i}$ \\
 & &  $u_{n+1}=(-1)^{n+1}\sqrt{(\alpha+1)w_{n+1}}$ \\
\hline
\end{tabular}
\end{center}

\caption{Relation of the barycentric weights with the Gauss (first row), Gauss--Radau (second and third rows) 
and Gauss--Lobatto (fourth row) weights.
 }
\end{table}

\section{Brief description of the MATLAB and Maple algorithms}

\begin{figure}[h]
\begin{framed}[\linewidth]
\begin{minipage}[t]{\linewidth}
\dirtree{%
.1 {\bf Matlab algorithms}.
.2 Gauss-Jacobi quadrature.
.3 \listtt{function GJ: Computation of Gauss-Jacobi quadrature}.
.3 \listtt{function legen: Computation of Gauss-Legendre quadrature}.
.3  \listtt{function rlJ: Computation of Gauss-Radau-Jacobi and Gauss-Lobatto-Jacobi quadratures and barycentric weights}.
.2 Gauss-Laguerre quadrature.
.3 \listtt{function GL: Computation of Gauss-Laguerre quadrature}.
.3  \listtt{function rL: Computation of Gauss-Radau-Laguerre quadrature and barycentric weights}.
.2 Gauss-Hermite quadrature.
.3 \listtt{function GHa: Asymptotic computation of Gauss-Hermite quadrature}.
.3  \listtt{function GHi: Iterative computation of Gauss-Hermite quadrature}.
.3  \listtt{function bH:  Computation of barycentric weights}.
.2 Examples.
.1 {\bf Maple algorithms}.
.2 Gauss-Gegenbauer quadrature.
.3 \texttt{\color{red}gegenbauer.mpl} \listtt{ Maple algorithm}.
.3 \texttt{\color{red}gegenbauer.mws}  \listtt{illustrative Maple worksheet}.
.2 Gauss-Hermite quadrature.
.3 \texttt{\color{red}hermite.mpl} \listtt{ Maple algorithm}.
.3 \texttt{\color{red}hermite.mws} \listtt{illustrative Maple worksheet}.
} 
\end{minipage}
\end{framed}
\caption{Overview of the software provided for Gauss quadratures}
\label{fig:software}
\end{figure}

The codes are available at

\url{http://github.com/NumericalQuadrature/GaussQuadrature}

\noindent
The original release (v1.0.0), publicly available at this site, contains the versions used in the present paper.

The different algorithms are documented at the website, including details on the input/output, modes of computation, 
and error handling; this information
is also provided in the comment lines in the headings of the codes. Additionally, 
examples showing the use of the MATLAB algorithms are provided in the MATLAB folder. Maple worksheets are included in the corresponding folder, demonstrating the use of the arbitrary-precision algorithms for the symmetric cases.

A description of the components of the software library is shown in Figure \ref{fig:software}.

\section*{Acknowledgments}
The authors gratefully acknowledge the two anonymous reviewers for their valuable comments and suggestions.

\bibliographystyle{siamplain}
\bibliography{gauss}
\end{document}